\documentclass{amsart}
\usepackage{amssymb,graphicx,upref,cite}
\usepackage[line,ps,color,poly,all,xdvi,dvips]{xy}
\newxyColor{lightgrey}{.75}{gray}{}
\CompileMatrices
\newcommand{\Z}{\mathbb Z}%
\DeclareMathOperator{\SU}{SU}%
\DeclareMathOperator{\Sym}{Sym}%
\DeclareMathOperator{\SL}{SL}%
\DeclareMathOperator{\Sp}{Sp}%
\DeclareMathOperator{\GL}{GL}%
\DeclareMathOperator{\vcd}{vcd}%
\newcommand{\R}{\mathbb R}%
\newcommand{\Q}{\mathbb Q}%
\newcommand{\C}{\mathbb C}%
\newcommand{\E}{\mathbb E}%
\newcommand{\G}{{\mathcal G}}%
\newcommand{\lsp}[1]{{}^{#1}\!}%
\DeclareMathOperator{\sgn}{sgn}%
\DeclareMathOperator{\rank}{rank}%
\DeclareMathOperator{\stab}{Stab}%
\newcommand{\Stab}[1]{\stab_{\Gamma}\left(#1\right)}%
\newcommand{\I} {{\mathcal I}}%
\newcommand{\F} {{\mathcal F}}%
\newcommand{\N} {{\mathcal N}}%
\newcommand{\Par} {{\mathcal P}}%
\newcommand{\QQ}[2]{{\mathcal Q}\left(#1,#2\right)}%
\newcommand{\QQQ}[1]{{\mathfrak Q}\left(\left[#1\right]\right)}%
\newcommand{\J} {{\mathcal J}}%
\renewcommand{\Re}{\operatorname{Re}}%
\renewcommand{\Im}{\operatorname{Im}}%
\newcommand{\conj}[1]{\lsp{#1}P_0}%
\newcommand{\vect}[3]{\begin{pmatrix} #1 \\ #2 \\ #3 \end{pmatrix}}%
\newcommand{\mat}[1]{\begin{pmatrix} #1 \end{pmatrix}}%
\newcommand{\ve}{\vect{1}{0}{0}}%
\newcommand{\vw}{\vect{0}{0}{1}}%
\newcommand{\vsigma}{\vect{i}{1+i}{1}}%
\newcommand{\comment}[1]{}%
\newcommand{\temp}[1]{}%
\newcommand{\Label}[1]{{\label{#1} \comment{#1}}}%
\newcommand{\vP}{(n,p,q)^t}%
\theoremstyle{plain}
\newtheorem{thm}{Theorem}[section]

\newtheorem{lem}[thm]{Lemma}
\newtheorem{prop}[thm]{Proposition}
\newtheorem{cor}[thm]{Corollary}
\theoremstyle{definition}
\newtheorem{defn}[thm]{Definition}

\begin{document}


\title{Explicit reduction for $\SU(2,1;\Z[i])$}
\author{Dan Yasaki}
\address{Department of Mathematics and Statistics\\Lederle Graduate Research Tower\\ University of Massachusetts\\Amherst, MA 01003-9305}
\email{yasaki@math.umass.edu}
\date{}
\thanks{The original manuscript was prepared with the \AmS-\LaTeX\ macro
system and the \Xy-pic\ package.}
\begin{abstract}
Let $\Gamma \backslash D$ be an arithmetic quotient of a symmetric space of non-compact type.  A spine $D_0$ is a $\Gamma$-equivariant deformation retraction of $D$ with dimension equal to the virtual cohomological dimension of $\Gamma$.  We explicitly construct a spine for the case of $\Gamma=\SU(2,1;\Z[i])$.  The spine is then used to compute the cohomology of $\Gamma\backslash D$ with various local coefficients.     
\end{abstract}
\maketitle

\bibliographystyle{../amsplain_initials}

\begin{section}{Introduction}\Label{sec:introduction}
Let $G$ be the real points of the $\Q$-rank 1 linear algebraic group $\SU(2,1)$, and let $D$ be the associated non-compact symmetric space.  Let $\Gamma$ be an arithmetic subgroup of the rational points $G(\Q)$.  Let $(E,\rho)$ be a $\Gamma$-module over $R$.  If $\Gamma$ is torsion-free, the locally symmetric space $\Gamma \backslash D$ is a $K(\Gamma,1)$ since $D$ is contractible, and the group cohomology of $\Gamma$ is isomorphic to the cohomology of the locally symmetric space, i.e. $H^*(\Gamma,E)\cong H^*(\Gamma \backslash D; \E)$, where $\E$ denotes the local system defined by $(E,\rho)$ on $\Gamma \backslash D$.  When $\Gamma$ has torsion, the correct treatment involves the language of orbifolds, but the isomorphism of cohomology is still valid by using a suitable sheaf $\E$ as long as the orders of the torsion elements of $\Gamma$ are invertible in $R$ 

The \emph{virtual cohomological dimension} ($\vcd$) of $G$ is the smallest integer $p$ such that cohomology of $\Gamma \backslash D$ vanishes in degrees above $p$, where $\Gamma\subset G(\Q)$ is any torsion-free arithmetic subgroup.  Borel and Serre \cite{BS} show that the discrepancy between the dimension of $D$ and the $\vcd(G)$ is given by the $\Q$-rank of $G$, the dimension of a maximal $\Q$-split torus in $G$.  Thus in our case, $D$ is 4-dimensional, and the virtual cohomological dimension of $\Gamma$ is $3$.  There is in fact a $3$-dimensional deformation retract $D_0\subset D$ that is invariant under the action of $\Gamma$ \cite{Yasrank1}.  Such spaces are known as \emph{spines}. 

Spines have been constructed for many groups \cite{Sou,Men,Vog,Brown,A,MM,Bat,LSz}.  In \cite{A2}, Ash describes the \emph{well-rounded retract}, a method for constructing a spine for all linear symmetric spaces.  Ash and McConnell \cite{AM} extend \cite{A2} to the Borel-Serre compactification and relate the retraction to a combination of geodesic actions.  The well-rounded retract has been used in the computation of cohomology \cite{Sou,LSz,A80,Men,Vog,SVog,AGG,AM1,AM2,vGT}.

The well-rounded retract proves the existence and gives a method of explicitly defining spines in linear symmetric spaces.  There were no non-linear examples until MacPherson and McConnell \cite{MM} constructed a spine in the Siegel upper half-space for the group $\Sp_4(\Z)$.  

In this paper, we provide another non-linear example by using the method of \cite{Yasrank1} to compute a spine for $\SU(2,1;\Z[i])$.  Sections~\ref{sec:su21background} and \ref{sec:construction} set notation and define the exhaustion functions that are used to describe the pieces of the spine.  In Section~\ref{sec:configuration}, we classify certain of configurations of isotropic line in $\C^3$.  We show the spine has the structure of a cell complex with cells related to these configurations in  Section~\ref{sec:admissible}.  Explicit $\Gamma$-representatives of cells are fixed, and their stabilizers are computed in Section~\ref{sec:stabilizers}.  After subdivision, we obtain a regular cell complex for $D_0$ on which $\Gamma$ acts cellularly.  In Section~\ref{sec:cohomology}, we recall some facts about orbifolds and develop machinery to investigate the cohomology of $\Gamma$.  The results of Section~\ref{sec:cohomology} hold in more generality, and may be of independent interest.  We apply these methods in Section~\ref{sec:cohcomputation} to $\Gamma=\SU(2,1;\Z[i])$ to compute the cohomology of $\Gamma$ with coefficients in various $\Gamma$-modules.

I would like to thank my thesis advisor, Les Saper for his insight into this work.  I would also like to thank Paul Gunnells for helpful conversations. 
\end{section}
\begin{section}{Preliminaries}\Label{sec:su21background}
Let $G$ be the identity component of the real points of the algebraic group $G=\SU(2,1)$, realized explicitly as
\[\SU(2,1)=\left\{g \in \SL(3,\C)\;\left |\; g^* \begin{pmatrix} 0&0&i\\0&-1&0\\-i&0&0 \end{pmatrix} g 
= \begin{pmatrix} 0&0&i\\0&-1&0\\-i&0&0 \end{pmatrix}\right.\right\}.\]  Alternatively, let $\mathcal{Q}$ be the (2,1)-quadratic form on $\C^3$ defined by 
\begin{equation*}
\QQ{u}{v}=u^*\begin{pmatrix} 0&0&i\\0&-1&0\\-i&0&0 \end{pmatrix} v.
\end{equation*}  
Then $G$ is the group of determinant 1 complex linear transformations of $\C^3$ that preserve $\mathcal Q$.  Let $\Gamma$ be the arithmetic subgroup $\Gamma = \SU(2,1)\cap \SL_3(\Z[i])$.

Let $\theta$ denote the Cartan involution given by inverse conjugate transpose and let $K$ be the fixed points under $\theta$.   Let $D=G/K$ be the associated Riemannian symmetric space of non-compact type.  Let $\Par$ denote the set of (proper) rational parabolic subgroups of $G$.

Let $P_0 \subset G$ be the rational parabolic subgroup of upper triangular matrices, and fix subgroups $N_0$, $A_0$, and $M_0$:
\begin{align*}
P_0&=\left\{\left.\begin{pmatrix} y\zeta & \beta \zeta^{-2} & \zeta\left(r + i|\beta|^2/2\right)/y \\ 0 & \zeta^{-2} & i{\overline \beta}\zeta/y \\ 0 & 0 & \zeta/y \end{pmatrix} \; \right| \;  \zeta,\beta \in \C,\ |\zeta|=1,\ r \in \R,\ y\in \R_{>0} \right\},\\
N_0&=\left\{\left.\begin{pmatrix} 1 & \beta & r + i|\beta|^2/2 \\ 0 & 1 & i{\overline \beta} \\ 0 & 0 & 1 \end{pmatrix} \; \right| \;  \beta \in \C,\ r \in \R\right\},\\
A_0&=\left\{\left.\begin{pmatrix} y&0&0\\0&1&0\\0&0&1/y \end{pmatrix}\; \right|\;  y \in \R_{>0} \right\},\\
M_0&=\left\{\left.\begin{pmatrix} \zeta&0&0\\0&\zeta^{-2}&0\\0&0&\zeta \end{pmatrix}\; \right|\;  \zeta \in \C, \ |\zeta|=1 \right\}.
\end{align*}
$P_0$ acts transitively on $D$, and every point $z\in D$ can be written as $p \cdot x_0$ for some $p \in P_0$.  Using  Langlands decomposition, there exists $u \in N_0, a \in A_0$, and $m \in M_0$ such that $p=uam$.  Since $M_0 \subset K$, $z$ can be written as $ua \cdot x_0$.  Denote such a point $z=(y,\beta,r)$.

Zink showed that $\Gamma$ has class number 1 \cite{Z}.  Thus $\Gamma \backslash G(\Q) / P_0(\Q)$ consists of a single point, and all the parabolic subgroups of $G$ are $\Gamma$-conjugate. The rational parabolic subgroups of $G$ are parametrized by the maximal isotropic subspaces of $\C^3$ which they stabilize.  These are $1$-dimensional, and so to each $P \in \Par$, there is an associated reduced, isotropic vector $v_P \in \Z[i]^3$.  (A vector $(n,p,q)^t \in \Z[i]^3$ is \emph{reduced} if $(n,p,q)$ generate $\Z[i]$ as an ideal.)  Similarly, given a reduced, isotropic vector $v$ in $\Z[i]^3$, there is an associated rational parabolic subgroup $P_v$.  Notice, however, that $v_P$ is only well-defined up to scaling by $\Z[i]^*=\{\pm 1, \pm i\}.$  Thus, the vectors $v$ and $\varepsilon v$ will be treated interchangeably for $\varepsilon \in \Z[i]^*$.  If $P=\lsp{\gamma}Q$ for some $\gamma \in \Gamma$, then $v_{P}=\gamma v_{Q}$.  

Unless explicitly mentioned otherwise, the vector $v_P$ will be written as $v_P=(n,p,q)^t$.  The isotropic condition $ \QQ{v_P}{v_P}=0$ implies that 
\begin{equation}\Label{eq:isotropic}
|p|^2=2\Im(n \overline{q}).
\end{equation}  
In particular, $q\neq 0$ for $P \neq P_0$.  Furthermore, since there are no isotropic 2-planes in $\C^3$,  
\begin{equation}
\QQ{v_P}{v_Q}\neq 0 \quad \text{for $P \neq Q$.}
\end{equation} 

Because these elements of $\Gamma$ will be used frequently, set once and for all\[w=\begin{pmatrix} 0 & 0 & -1 \\ 0 & 1 & 0 \\ 1 & 0 & 0 \end{pmatrix},
\sigma=\begin{pmatrix} 1 & 1+i & i \\ 0 & 1 & 1+i \\ 0 & 0 & 1 \end{pmatrix}, \check{\sigma}=\begin{pmatrix} 1 & i(1+i) & i \\ 0 & 1 & -i(1+i) \\ 0 & 0 & 1\end{pmatrix},\] \[\tau= \begin{pmatrix} 1 & 0 & 1 \\ 0 & 1 & 0 \\ 0 & 0 & 1 \end{pmatrix}, \epsilon=\begin{pmatrix} i & 0 & 0 \\ 0 & -1 & 0 \\ 0 & 0 & i \end{pmatrix},\text{ and $\xi=\tau w \tau \sigma w \epsilon^3$.}\]  Note that $\check{\sigma}$ is contained in the group generated by $\{\epsilon, w, \sigma\}$.  In particular, $\check{\sigma}=w \sigma\epsilon w \sigma^{-1}w$.
\end{section}
\begin{section}{Construction of the spine}\Label{sec:construction}
In this section we briefly describe the construction of a $\Gamma$-invariant, 3-dimensional cell complex which is a deformation retract of $D$.  This construction is described for the general $\Q$-rank 1 case in \cite{Yasrank1}.  

We first define an exhaustion function $f_P$ for every rational parabolic subgroup $P \subseteq G$.  These exhaustion functions are then used to define a decomposition of $D$ into sets $D(\I)$ for $\I \subset \Par$. 
\begin{subsection}{Exhaustion functions}
Let $z=(y,\beta,r)\in D$ and $P$ a rational parabolic subgroup of $G$ with associated isotropic vector $\vP$. Then the exhaustion function $f_P$ can be written as
\begin{align}\Label{eq:an}
 f_{0}(z)&\equiv  f_{P_0}(z)=y\\
 f_P(z)&=\frac{y}{\left(|n-\beta p+\bigl(i|\beta|^2/2-r\bigr)q|^2+y^2|p-i\overline \beta q|^2+y^4|q|^2\right)^{1/2}}.
\end{align}
The family of exhaustion functions defined above is $\Gamma$-invariant in the sense that 
\begin{equation}\Label{eq:invariantfamily}
f_{\lsp{\gamma}P}(z)=f_P(\gamma^{-1}\cdot z)\quad \text{for $\gamma \in \Gamma$.}
\end{equation}

\end{subsection}
\begin{subsection}{Admissible sets}
For a parabolic $P$, define $D(P)\subset D$ to be the set of $z \in D$ such that $ f_P(z) \geq  f_Q(z)$ for every $Q \in \Par \setminus \{P\}$.  More generally, for a subset  $\I \subseteq \Par$, 
\begin{align}
E(\I)&=\{z\in D\; |\;  f_P(z)= f_Q(z) \text{ for every pair } P,Q \ \in \I \}\\
D(\I)&=\bigcap_{P\in \I} D(P)\\ 
D'(\I)&=D(\I) \setminus \bigcup_{\I' \supsetneq \I} D(\I').
\end{align}
It follows that $D'(\I)\subseteq D(\I) \subset E(\I) \text{ and } D(\I)=\coprod_{\tilde{\I} \supseteq \I}D'(\tilde{\I})$.  Let $f_\I$ denote the restriction to $E(\I)$ of $f_P$ for $P\in \I$. 
\begin{defn}
Let $\I \subseteq \Par$, $P \in \I$, and $z \in E(\I)$.  Then $z$ is called a {\it first contact for $\I$} if $f_\I(z)$ is a global maximum of $f_\I$ on $E(\I)$.
\end{defn}
\begin{defn}
A subset $\I \subset \Par$ is called {\it admissible} if $D(\I)$ is non-empty and {\it strongly admissible} if $D'(\I)$ is non-empty.
\end{defn}

Let $D_0 \subset D$ denote the subset
\begin{equation}\Label{eq:spine}
D_0=\coprod_{|\I| >1} D'(\I).
\end{equation} 
The deformation is defined separately on each $D(P)$ for $P\in \Par$.  For $z\in D(P)$, we use the (negative) gradient flow of $f_P$ to flow $z$ to a point on $D_0$.  This corresponds to using the geodesic action \cite{BS} of $A_P$ on $z$ \cite{Yasrank1}.

Let $f_{D_0}$ denote the function on $D_0$ given by
\begin{equation}\Label{eq:fD0}
f_{D_0}(z)=f_\I(z) \quad \text{for $z\in D(\I)$.}
\end{equation}

\end{subsection}
\begin{subsection}{First contact points}
Given the explicit description of the exhaustion functions in coordinates, one can readily describe the set $E(\{P_0,P\})$.  Writing $z=(y,\beta,r)$ and using \eqref{eq:an},
\begin{equation} \Label{eq:phi}
\left(\frac{ f_0(z)}{ f_P(z)}\right)^2= \left|n-\beta p+\left(i|\beta|^2/2-r\right)q\right|^2+y^2\left|p-i\overline \beta q\right|^2+y^4|q|^2.
\end{equation}

\begin{prop}\Label{prop:E(P0,P)}
Let $P$ be a rational parabolic subgroup of $G$ with associated isotropic vector $\vP$.  Every $z=(y,\beta,r)\in E(\{P_0,P\})$ satisfies
\begin{equation*}
y^2=-\frac{1}{2}\left|\frac{p}{q}-i\bar \beta\right|^2+\sqrt{\frac{1}{|q|^2}-\left(\Re\left(\frac{n-\beta p}{q}-r\right)\right)^2}
\end{equation*}
\end{prop}
\begin{proof}
Note that for $z \in E(\{P_0,P\})$, $\left(\frac{ f_0(z)}{ f_P(z)}\right)^2=1$.  Then the result follows from \eqref{eq:phi} using the quadratic formula to solve for $y^2$, and simplifying the result using the isotropic condition \eqref{eq:isotropic}.
\end{proof}

Proposition~\ref{prop:E(P0,P)} allows us to easily calculate the first contact point for $\{P_0,P\}$.  The $\Gamma$-invariance of the exhaustion functions allows us to translate this for general $\{P,Q\}$. 
\begin{prop} \Label{prop:first}
Let $\I=\{P,Q\} \subset \Par$.  Let $z$ be a first contact point for $\I$.  Then 
\begin{equation*}
 f_P(z)= f_Q(z)=\frac{1}{\sqrt{|\QQ{v_P}{v_Q}|}}.
\end{equation*}
In particular, the first contact for $\{P_0,P\}$ is $z=\left(1/\sqrt{|q|},i(\overline{p/q}),\Re(n/q)\right)$.  
\end{prop}
\end{subsection}
\end{section}
\begin{section}{Configurations of vectors}\Label{sec:configuration}
\begin{defn}
Let $\J$ be a subset of vectors in $\C^3$.  Then $\J$ is said to be {\it $c$-bounded} if 
\[|\QQ{u}{v}|^2\leq c \quad \text{for every $u$ and $v$ in $\J$.}\]
\end{defn}

Set once and for all the following sets $\J^i_j$ of isotropic vectors in $\C^3$.
\begin{gather*}
\J^2_1=\left\{\vect{1}{0}{0},\vect{0}{0}{1}\right\} \qquad\J^2_2=\left\{\vect{1}{0}{0},\vect{i}{1+i}{1+i}\right\}\\
\J^3_1=\J^2_1 \cup \left\{\vect{1}{0}{1}\right\}\qquad\J^3_2=\J^2_1 \cup\left\{\vect{i}{1+i}{1}\right\}\qquad \J^3_3=\J^2_1  \cup\left\{\vect{1+i}{1+i}{1}\right\}\\
\J^4_1=\J^3_1\cup\left\{\vsigma \right\} \qquad \J^4_2=\J^3_3\cup \left\{\vect{-1}{-1+i}{1+i}\right\}\\
\J^5=\J^4_1 \cup\left\{\vect{1+i}{1+i}{1}\right\}\\
\J^8=\left\{\ve,\vw,\begin{pmatrix} -1\\1+i\\1+i \end{pmatrix},\vect{-1+i}{1+i}{1},\vect{1+i}{1-i}{1},\begin{pmatrix} i\\1+i\\1+i \end{pmatrix},\vect{2i}{2}{1},\vect{i}{2}{2}\right\}.
\end{gather*}
One easily checks that these sets are not $\Gamma$-conjugate, and $\J^8$ is $4$-bounded, while the other sets are $2$-bounded.
\begin{prop}\Label{prop:configurations}
A $2$-bounded set of reduced, integral, isotropic vectors is $\Gamma$-conjugate to exactly one of $\J^2_1,\ \J^2_2,\ \J^3_1,\ \J^3_2,\ \J^3_3,\ \J^4_1,\ \J^4_2,$ or $\J^5$.
\end{prop}

\begin{proof}
Let $\J$ denote such a subset of order 2.  Since $G$ has class number $1$, we can assume that one of the vectors of $\J$ is $(1,0,0)^t$.   Let $v=\vP$ be the other vector of $\J$.  Let $\sigma, \ \check{\sigma}, \ \tau , \text{ and } \epsilon$ be defined as in Section~\ref{sec:su21background}.  These elements preserve $(1,0,0)^t$ (up to scaling by $\Z[i]^*$).  By applying powers of $\sigma$ and $\check{\sigma}$ to $v$, we can add any $\Z$-linear combination of $(1+i)q$ and $(1+i)iq$ to $p$ to force it into the square in the complex plane with vertices $q, \ iq, \ -q,$ and $-iq$.  By applying powers of $-i\epsilon$ to $v$, one can force $p$ to lie in the triangle with $q, \ iq,$ and 0 as vertices while leaving $q$ fixed.  Then by applying powers of $\tau$ to $v$, one can add a $\Z$-scalar multiple of $q$ to $n$ so that $n$ now has the form $dq+ciq$, where $-\frac{1}{2}<d\leq \frac{1}{2}$ for some $c \in \R$.  Since $v$ is isotropic, $|p|^2=2\Im(n \overline q)=2c|q|^2$.  Since $\J$ is $2$-bounded, $|q|^2\leq 2$, so that in particular, $\J$ is $\Gamma$-equivalent to $\J^2_1$ or $\J^2_2$.  

For $|\J|>2$, we can arrange that  $\J \supset \J^2_1$ or $\J \supset \J^2_2$.  The $2$-bounded condition in each of these two cases has only finitely many solutions.  The proposition then follows from listing the solutions and identifying $\Gamma$-conjugate sets.
\end{proof}
\end{section}
\begin{section}{Reduction theory}\Label{sec:reduction}
\begin{prop}\Label{prop:Siegelstrip}
Every point $z \in D$ is conjugate under $\Gamma_{P_0}$ to a point $(y,\beta,r)$, where $-\frac{1}{2}<r \leq \frac{1}{2}$ and $\beta$ lies in the square in the complex plane with vertices  $0, \ \frac{1+i}{2}, \ i,$ and $\frac{-1+i}{2}$. 
\end{prop}
\begin{proof}
Consider the elements $\{\sigma,~\check{\sigma},~\tau,~\epsilon\} \subset \Gamma_{P_0}$ defined in Section~\ref{sec:su21background}.  The action of $\{\sigma,~\check{\sigma},~\tau,~\epsilon\}$ leaves $D(P_0)\cap D_0$ stable, and is given explicitly by
\begin{equation}
\begin{split}
\sigma \cdot (y,\beta,r)&=(y,\beta+(1+i),r-\Re(\beta)+\Im(\beta)),\\ 
\check{\sigma} \cdot (y,\beta,r)&=(y,\beta+i(1+i),r-\Re(\beta)-\Im(\beta)),\\
\tau \cdot (y,\beta,r)&=(y,\beta,r+1), \quad \text{and}\\
\epsilon \cdot (y,\beta,r)&=(y,-i\beta,r).   
\end{split}
\end{equation}   
Thus, by applying powers of $\sigma \ \text{and }\check{\sigma}$, $\beta$ can be put in the square in the complex plane with vertices $1,\ i,\ -1,\ \text{and } -i$.  By applying a power of $\epsilon,$ $\beta$ can be put in the square in the complex plane with vertices $0, \ \frac{1+i}{2}, \ i, \ \text{and } \frac{-1+i}{2}$.  Then by applying powers of $\tau$, it can be arranged that $-\frac{1}{2} < r \leq \frac{1}{2}$ without changing the value of $\beta$.
\end{proof}

\begin{prop}\Label{prop:mumin}
For each $z\in D_0$, 
\[\frac{1}{\root 4 \of 5}<f_{D_0}(z)\leq 1.\]
\end{prop}
\begin{proof}
Since $f_{D_0}(z)=f_{D_0}(\gamma \cdot z)$ for $z \in D_0$ and $\gamma \in \Gamma$, and every point in $D_0$ is $\Gamma$-conjugate to a point of $D(P_0)$, it suffices to determine the range of $f_{D_0}$ on $D(P_0)\cap D_0$.  In fact, it suffices to determine the range on a subset $F \subseteq D(P_0)\cap D_0$, provided the $\Gamma$-translates of $F$ cover $D(P_0)\cap D_0$.  The action of $\{\sigma, \check{\sigma}, \tau, \epsilon\}$ leaves $D(P_0)\cap D_0$ stable, so it suffices to determine the range on $F=D(P_0)\cap D_0 \cap T$, where $T$ is the strip in $D$ defined in Proposition~\ref{prop:Siegelstrip}.  By construction, $f_{D_0}(z)=\max_{P \in \Par}\{ f_P(z)\}$ and Proposition~\ref{prop:first} shows that $f_{D_0}(z)\leq 1$.  
Thus it suffices to show that 
\begin{equation}
\min_{z \in F}\{ f_0(z)\} > \frac{1}{\root 4 \of 5}.
\end{equation}
Note that on $F$, $ f_0(z)\geq  f_P(z)$ for all $P \in \Par$.  By Proposition~\ref{prop:E(P0,P)}, this is equivalent to the condition that 
\begin{equation}\Label{eq:ybig}
 f_0(z)^2=y^2 \geq -\frac{1}{2}\left |\frac{p}{q}-i \overline{\beta}\right |^2+\sqrt{\frac{1}{|q|^2}-\left(\Re\left(\frac{n-\beta p}{q}-r\right)\right)^2}.
\end{equation}

Consider the rational parabolic subgroups $\conj{w}$, $P$, and $Q$ corresponding to the vectors $v_w=(0,0,1)^t, \ v_P=(i, 1+i, 1+i)^t\ \text{and } v_Q=(-1, 1+i , 1+i)^t$.  Divide $F$ into the following regions:
\begin{align*}
F_P&=\Bigl\{x=(y,\beta,r) \in F\ \Bigl|\ |\beta|>\frac{9}{10},\ \frac{3}{10} <r\leq \frac{1}{2}\Bigr\}\\
F_Q&=\Bigl\{x=(y,\beta,r) \in F\ \Bigl|\ |\beta|>\frac{9}{10},\ -\frac{1}{2} \leq r<-\frac{3}{10}\Bigr\}\\
F_w&=F \setminus \{F_1 \cup F_2\}.
\end{align*}

Since $\beta$ is constrained to lie in the square indicated in Proposition~\ref{prop:Siegelstrip}, one may calculate that on $F_P$ and $F_Q$
\begin{gather}
-\frac{(10-\sqrt{62})}{20} <\Re(\beta)<\frac{10-\sqrt{62}}{20}\\
|\beta-i|^2<\frac{81-10\sqrt{62}}{100}.  
\end{gather}
Comparing $ f_0$ to $ f_{P}$ on $F_P$, $f_0$ to $f_{Q}$ on $F_Q$, and $ f_0$ to $ f_{w}$ on $F_w$, \eqref{eq:ybig} gives the desired bound on each piece.
Since $F=F_P \cup F_Q \cup F_w$, this proves the result.
\end{proof}
\end{section}

\begin{section}{Admissible sets}\Label{sec:admissible}
Since $v_P$ and $v_Q$ are integral vectors, $\QQ{v_P}{v_Q} \in \Z[i]$.  Proposition~\ref{prop:first} and Proposition~\ref{prop:mumin} imply the following.  
\begin{prop} \Label{prop:Qbound}
Let $\I$ be an admissible set.  Then 
\[|\QQ{v_P}{v_Q}|^2 \leq 4 \quad \text{for every $P$ and $Q$ in $\I$.}\]
In particular, the set of vectors associated to an admissible set is $4$-bounded.
\end{prop}
\begin{prop} \Label{prop:qsize}
Let $\I =\{P,Q\}$ be an admissible subset of $\Par$.  If $|\QQ{v_P}{v_Q}|^2=4$, then there exists a strongly admissible set $\tilde \I \subset \Par$ of order $8$ such that $D(\I)= D(\tilde \I)=\{z\}$, where $z$ is the first contact for $\I$.  Let $\tilde \J$ be the set of isotropic vectors associated to $\tilde \I$.  Then $\tilde \J$ is $4$-bounded and is $\Gamma$-equivalent to $\J^8$.  A set of isotropic vectors that is not $2$-bounded, but is associated to a strongly admissible set, is $\Gamma$-equivalent to $\J^8$. 
\end{prop}
\begin{proof}
Suppose  $\I =\{P,Q\}$ is admissible and $|\QQ{v_P}{v_Q}|^2=4$.  Since $G$ has class number 1, $\I$ is $\Gamma$-conjugate to a set of the form $\{P_0,P\}$.  By $\Gamma$-action,  one only needs to consider the cases when $v_P=(i,2,2)^t$ and $v_P=(1,0,2)^t$. 

In the case that $v=(1,0,2)^t$ we claim that $\{P_0,P_v\}$ is not admissible, and hence $D(\{P_0,P_v\})$ is empty.  To see this, it suffices to show that $ f_Q(x)> f_0(x)$ on $E(\{P_0,P_v\})$ for some $Q$.  Proposition~\ref{prop:E(P0,P)} implies that $E(\{P_0,P_v\})$ is the set where 
\begin{equation}\label{eq:ov}
y^2=-\frac{|\beta|^2}{2}+\sqrt{r-r^2}.
\end{equation}
In particular, $0<r<1$.  Explicit computation with \eqref{eq:an} and \eqref{eq:ov} shows that on $E(\{P_0,P_v\})$, for $Q=\lsp{w}P_0$, $f_Q(z)=\frac{ f_0(z)}{\sqrt{r}}> f_0(z)$.

Next, consider the case $v=(i,2,2)^t$.  Proposition~\ref{prop:E(P0,P)} implies that $E(\{P_0,P_v \})$ is the set where 
\begin{equation}\Label{eq:i22}
y^2=-\frac{|\beta-i|^2}{2}+\frac{1}{2}\sqrt{1-4(\Re(\beta)+r)^2}.
\end{equation}
Consider the rational parabolic subgroups $P_1$ and $P_2$ with associated isotropic vectors $( i,1+i,1+i)^t$ and $( -1,1+i,1+i)^t$ respectively.  Explicit computation with \eqref{eq:an} and \eqref{eq:i22} shows that on $E(\{P_0,P_v\})$,
 \[f_{P_1}(x)=\frac{ f_0(x)}{1-2r-\Re(\beta)} \quad \text{and} \quad f_{P_2}(x)=\frac{ f_0(x)}{1+2r+\Re(\beta)}.\]
Thus on $D(\{P_0,P_v\})$, where $f_{P_i}(z)\leq f_0(z)$, we see that $2r=-\Re(\beta)$, and so
\begin{equation*}
D(\I)=D(\{P_0,P_v,P_1,P_2\}).
\end{equation*}
In particular, note that $|\beta-i|^2<\sqrt{1-4r^2}\leq 1$.  

Consider the rational parabolic subgroups $Q_1,\ Q_2,\ Q_3,\ \text{and } Q_4$ with associated isotropic vectors $(0,0,1)^t$, $(1+i,1-i,1)^t$, $(2i,2,1)^t$, and $(-1+i,1+i,1)^t$
 respectively.  Explicit computation shows that on $S$,
\begin{align*}
 f_{Q_1}(x)&=\frac{ f_0(x)}{|\beta|^2}&
 f_{Q_2}(x)&=\frac{ f_0(x)}{|\beta-(-1+i)|^2}\\
 f_{Q_3}(x)&=\frac{ f_0(x)}{|\beta-2i|^2}
& f_{Q_4}(x)&=\frac{ f_0(x)}{|\beta-(1+i)|^2}.\\
\end{align*}
Since $|\beta-i|^2<1$ and $ f_{Q_j}(x)\leq  f_0(x)$ on $D(\I)$,  we see that $\beta$ is forced to equal $i$.  Thus 
\begin{equation*}
D(\I)=D(\{P_0,P_v,P_1,P_2,Q_1,Q_2,Q_3,Q_4\})=\left\{\left(\frac{1}{\sqrt{2}},i,0\right)\right\}.
\end{equation*}

One checks that $\J^8$ is not properly contained in any $4$-bounded set to complete the proof.
\end{proof}

\begin{cor}\Label{cor:lessthan8}
The set of isotropic vectors associated to a strongly admissible set of order less than eight is $2$-bounded.  Furthermore, there are no strongly admissible sets of order greater than eight.
\end{cor}
\begin{proof}
Propositions~\ref{prop:Qbound}~and~\ref{prop:qsize} imply the first statement.  The second follows from explicit calculation that shows there are no $4$-bounded sets that properly contain $\J^8$.
\end{proof}

\begin{prop}\Label{prop:I}
Let $\I^i_j$ denote the set of rational parabolic subgroups associated $\J^i_j$ as defined in Section~\ref{sec:configuration}. Then
\begin{gather*}
\I^2_1=\{P_0,\conj{w}\} \qquad \I^2_2=\{P_0,\conj{\xi}\}\\
\I^3_1=\I^2_1 \cup\{\conj{\tau w}\} \qquad \I^3_2=\I^2_1 \cup \{\conj{\sigma w}\} \qquad \I^3_3=\I^2_1 \cup \{\conj{\tau \sigma w}\} \\
\I^4_1=\I^3_1 \cup \{\conj{\sigma w}\}\qquad \I^4_2=\I^3_3 \cup\{\conj{w^{-1}\tau \check{\sigma} w}\}\\
\I^5=\I^4_1 \cup \{\conj{\tau \sigma w}\}\\
\I^8=\{P_0,\conj{w},\conj{w \tau \sigma w},\conj{\tau^{-1} \sigma w},\conj{\tau \check{\sigma} w},\conj{\tau w \tau \sigma w},\conj{\tau^2 \check{\sigma} \sigma w},\conj{\epsilon w \xi^4 w}\}
\end{gather*}
\end{prop}
Using \eqref{eq:an} and Propositions~\ref{prop:first} and \ref{prop:I}, one calculates the following first contact points $z(\I)$.
\begin{prop} \Label{prop:explicitfirst}
Let $\phi=\frac{1+\sqrt{5}}{2}$.  Then
\begin{gather*}
z(\I^2_1)=(1,0,0)\qquad z(\I^2_2)=\left(\frac{1}{\sqrt[4]{2}},i,\frac{1}{2}\right)\\
z(\I^3_1)=\biggl(\sqrt[4]{\frac{3}{4}},0,\frac{1}{2}\biggr)\qquad z(\I^3_2)=\biggl(\frac{\sqrt{3}}{2},\frac{1+i}{2},0\biggr)\\
z(\I^3_3)=\left(\frac{1}{\sqrt{2}}\sqrt{\phi^2\sqrt{\phi}-2},\frac{1}{2}(\phi^2-\sqrt{\phi}+i(1-\phi+\sqrt{\phi})),\frac{1}{2}(1-\phi+\sqrt{\phi})\right)\\
z(\I^4_1)=\biggl(\frac{\sqrt{-3+\sqrt{3}+\sqrt{2}+\sqrt{6}}}{2},\frac{1+\sqrt{3}-\sqrt{2}}{4}(1+\sqrt{3}i),\frac{1}{2}\biggr)\\
z(\I^4_2)=\Bigl(\sqrt{\phi-1},\frac{3+i}{2}(1-\sqrt{2}),0\Bigr)\\
z(\I^5)=\biggl(\frac{\sqrt{-1+2\sqrt{3}}}{2},\frac{1+i}{2},\frac{1}{2}\biggr)\qquad z(\I^8)=\left(\frac{1}{\sqrt{2}},i,0\right).
\end{gather*}
\end{prop}
\begin{thm}
Up to $\Gamma$-conjugacy, the strongly admissible sets are exactly \[\I^2_1,\ \I^2_2,\ \I^3_1,\ \I^3_2,\ \I^3_3,\ \I^4_1, \ \I^4_2,\ \I^5,\ \I^8.\]
\end{thm}
\begin{proof}
Propositions~\ref{prop:configurations}~and~\ref{prop:qsize} and Corollary~\ref{cor:lessthan8} imply that every strongly admissible set must be one of the ones listed.  Thus, it suffices to show that for each set $\I$ listed, $D'(\I)$ is non-empty.  In particular, it suffices to show that $D'(\I)$ contains its first contact point.  One uses \eqref{eq:an} and Proposition~\ref{prop:explicitfirst} to show that for $P\in \Par \setminus \I$, $f_P(z(\I))<f_\I(z(\I))$.  For example, $z(\I^2_1)=(1,0,0)$ and   
\[f_P(1,0,0)=\frac{1}{\sqrt{|n|^2+|p|^2+|q|^2}}<1 \quad \text{for all $P\in \Par \setminus \I^2_1$}.\]
The other cases follow similarly.
\end{proof}
\end{section}

\begin{section}{Pieces of the spine}\Label{sec:pieces}
\begin{thm}\Label{thm:pieces} Let $(y,\beta,r)$ denote a point in $D$.  Then 
\begin{enumerate}
\item $E(\I_1^2)=\left\{y^2=-\frac{|\beta|^2}{2}+\sqrt{1-r^2}\right\}$\Label{it:1}
\item $E(\I^2_2)=\left\{y^2=-\frac{1}{2}\left|1-i\bar{\beta}\right|^2+\sqrt{\frac{1}{2}-\left(\frac{1}{2}-\Re(\beta)-r\right)^2}\right\}$\Label{it:2}
\item $D(\I_1^2)\subseteq E(\I^2_1)\cap \left\{\frac{1}{\sqrt{5}} < y^2 \leq 1,\ -\frac{1}{2}<r\leq\frac{1}{2},\ |\beta|^2<-\frac{2}{\sqrt{5}}+2 \right\}$\Label{it:I21}
\item $D(\I^2_2) \subseteq E(\I^2_2) \cap \left\{0\leq \Re(\beta)+r <\sqrt{\frac{3}{10}} +\frac{1}{2}\right\}$ \Label{it:I22}
\item $D(\I^3_1)=D(\I^2_1)\cap \left\{r=\frac{1}{2}\right\}$\Label{it:I31}
\item $D(\I^3_2)=D(\I^2_1)\cap \left\{y^2=-\frac{1}{2}|1+i-i\bar{\beta}|^2+\sqrt{1-(\Re(\beta(1+i))+r)^2}\right\}$\Label{it:I32}
\item $D(\I^3_3)=D(\I^2_1)\cap \left\{y^2=-\frac{1}{2}|1+i-i\bar{\beta}|^2+\sqrt{1-(1-r-\Re(\beta(1+i)))^2}\right\}$\Label{it:I33}
\item $D(\I^4_1)=D(\I^2_1)\cap\left\{r=\frac{1}{2},\,\left|\beta-\frac{1+\sqrt{3}}{4}(1+\sqrt{3}i)\right|^2=\frac{1}{2}\right\}$\Label{it:I41}
\item $D(\I^4_2)=D(\I^2_1)\cap \left\{r=0,\, \left|\beta-\frac{3+i}{2}\right|^2=\frac{1}{2}\right\}$\Label{it:I42}
\item $D(\I^5)=\left\{\left(\frac{\sqrt{-1+2\sqrt{3}}}{2},\frac{1+i}{2},\frac{1}{2}\right)\right\}$\Label{it:I5}
\item $D(\I^8)=\left\{\left(\frac{1}{\sqrt{2}},i,0\right)\right\}$\Label{it:I8}
\end{enumerate}
\end{thm}
\begin{proof} 
Proposition~\ref{prop:E(P0,P)} implies \eqref{it:1}~and~\eqref{it:2}.  Similarly, \eqref{it:I31}-\eqref{it:I8} follow from computer calculations and repeated uses of Proposition~\ref{prop:E(P0,P)}.  

To show the bounds in \eqref{it:I21}, let $P=\conj{w}$ and consider the rational parabolic subgroups $Q=\conj{\tau}$ and $R=\conj{\tau^{-1}}$ with associated isotropic vector $v_Q=(1,0,1)^t$ and $v_R=(-1,0,1)^t$ respectively.  Then \eqref{eq:an} implies that
\begin{align}\Label{eq:f1}
f_{P}(z)&=\frac{y}{\left(|i\beta|^2/2-r|^2+y^2|\beta|^2+y^4\right)^{1/2}}\\\Label{eq:f2}
f_Q(z)&=\frac{y}{\left(|1 +i|\beta|^2/2-r|^2+y^2|\beta|^2+y^4\right)^{1/2}}\\\Label{eq:f3}
f_R(z)&=\frac{y}{\left(|-1 +i|\beta|^2/2-r|^2+y^2|\beta|^2+y^4\right)^{1/2}}
\end{align} 
For $z\in D(\I^2_1)$, $f_{P}(z) \geq f_Q(z)$ and $f_{\conj{w}}(z) \geq f_Q(z)$ and hence \eqref{eq:f1}, \eqref{eq:f2}, and \eqref{eq:f3} imply that $-\frac{1}{2} \leq r \leq \frac{1}{2}$. Since $f_0(z)=y$, Proposition~\ref{prop:mumin} implies that $\frac{1}{\sqrt{5}} <y^2 \leq 1$.  It follows that $|\beta|^2< -\frac{2}{\sqrt{5}}+2$ on $D(\I^2_1)$.  

To show \eqref{it:I22}, note that $y^2>\frac{1}{\sqrt{5}}$ on $D(\I^2_2)$ by Proposition~\ref{prop:mumin}, and hence $\Re(\beta)+r < \sqrt{\frac{3}{10}}+\frac{1}{2}$.  Let $P=\conj{\xi}$ and consider the rational parabolic subgroup $Q$ associated to the isotropic vector $v_Q=(-1,1+i,1+i)^t$.  For $z\in D(\I^2_2)$, $f_P(z)=y$ and $f_P(z) \geq f_Q(z)$ and hence \eqref{eq:an} implies that 
\begin{equation}\Label{eq:stuff}
|i-\beta(1+i)+(i|\beta|^2/2-r)(1+i)|^2\leq |-1-\beta(1+i)+(i|\beta|^2/2-r)(1+i)|^2.
\end{equation}
Thus \eqref{eq:stuff} and \eqref{eq:an} imply $0 \leq \Re(\beta)+r$.
\end{proof}

From the explicit description of representatives of $\Gamma$-conjugacy classes of admissible sets given above and the exhaustion functions given in \eqref{eq:an}, one can calculate a strict lower bound of $f_{D_0}$.
\begin{prop}\Label{prop:strongmumin}
For each $z\in D_0$,
\[
\frac{1}{\sqrt{2}}\leq f_{D_0}(z) \leq 1.
\]
\end{prop}
\end{section}

\begin{section}{Stabilizers}\Label{sec:stabilizers}
\begin{prop}\Label{prop:I21stab}
The stabilizer in $\Gamma$ of $D'(\I_1^2)$ is $\{e,\epsilon,\epsilon^2,\epsilon^3,w,\epsilon w,\epsilon^2 w,\epsilon^3 w\}$ and is isomorphic to $\Z/4\Z \times \Z/2\Z$ via the morphism which sends $\epsilon$ to the generator of $\Z/4\Z$ and $\epsilon w$ to the generator of $\Z/2\Z$.  
\end{prop}

\begin{proof}
Note that $\Stab{D'(\I_1^2)}=\Stab{\I^2_1}$.  One can calculate that $w^2\in \Gamma_{P_0}$ and hence $w \in \Stab{\I^2_1}$ and acts non-trivially on the ordered pair $(P_0,\conj{w})$.  Thus, $\Stab{\I^2_1}=L \cup L\cdot w$, where $\gamma \in L$ if and only if $\conj{\gamma}=P_0$ and $\conj{\gamma w}=\conj{w}$.  Since a parabolic subgroup is its own normalizer, one computes that 
\[
L= \Gamma \cap P_0 \cap \conj{w}=\{e,\epsilon,\epsilon^2,\epsilon^3\},
\]
and hence, 
\[
\Stab{D'(\I^2_1)}= \{e,\epsilon,\epsilon^2,\epsilon^3,w,\epsilon w,\epsilon^2 w,\epsilon^3 w\}.
\]

  It is easily checked that the map given in the proposition is an isomorphism.
\end{proof}

\begin{prop}\Label{prop:I22stab}
The stabilizer in $\Gamma$ of $D'(\I_2^2)$ is the cyclic group of order eight generated by $\xi$.
\end{prop}
\begin{proof}
Note that $\Stab{D'(\I_2^2)}=\Stab{\I^2_2}$.  Recall that $\I^2_2=\{P_0,\conj{\xi}\}$.  One can calculate that $\xi^2=\epsilon \tau \sigma^{-1} \in \Gamma_{P_0}$ and hence $\xi$ is in $\Stab{\I^2_2}$ and acts non-trivially on the ordered pair $(P_0,\conj{\xi})$.  Thus, $\Stab{\I^2_2}=L\cup L\cdot \xi$, where $\gamma \in L$ if and only if $\conj{\gamma}=P_0$ and $\conj{\gamma \xi}=\conj{\xi}$.  Since a parabolic subgroup is its own normalizer, one computes that
\begin{equation*}
L= \Gamma \cap P_0 \cap \conj{\xi}=\{e,\xi^2,\xi^4,\xi^6\},
\end{equation*} 
and hence $\Stab{D'(\I^2_2)}$ is the cyclic group of order eight generated by $\xi$. \end{proof}

\begin{prop}\Label{prop:I31stab}
The stabilizer in $\Gamma$ of $D'(\I^3_1)$ is the cyclic group of order twelve generated by $\tau \epsilon w$.
\end{prop}
\begin{proof}
Note that $\Stab{D'(\I^3_1)}=\Stab{\I^3_1}$.  One can easily check that 
\[
{}^\tau\{P_0,\conj{w}\}=\{P_0,\conj{\tau w}\} \quad \text{and}\quad{}^{w \tau^{-1} w}\{P_0,\conj{w}\}=\{\conj{w},\conj{\tau w}\}.
\]
Thus $\Stab{\I^3_1}=\Gamma_1 \cup \Gamma_2 \cup \Gamma_3$,
where
\begin{align*}
\Gamma_1&=\left\{\left.\gamma \in \Stab{\I^2_1}\;\right|\; \conj{\gamma \tau w}=\conj{\tau w}\right\},\\
\Gamma_2&=\left\{\left.\gamma \in \tau \cdot \Stab{\I^2_1}\;\right|\;\conj{\gamma \tau w}=\conj{w}\right\}, \quad \text{and}\\
\Gamma_3&=\left\{\left.\gamma \in w \tau^{-1} w \cdot \Stab{\I^2_1}\;\right|\; \conj{\gamma \tau w}=P_0\right\}.
\end{align*}
One calculates using Proposition~\ref{prop:I21stab} that 
\begin{align*}
\Gamma_1&=\{e,\epsilon,\epsilon^2,\epsilon^3\},\\
\Gamma_2&=\{\tau w, \tau \epsilon w,\tau \epsilon^2 w,\tau \epsilon^3 w\}, \quad \text{and}\\
\Gamma_3&=\{w \tau^{-1},w\tau^{-1}\epsilon,w\tau^{-1}\epsilon^2,w\tau^{-1}\epsilon^3\}.
\end{align*}
Explicit matrix multiplication shows that $\Gamma_1 \cup \Gamma_2 \cup \Gamma_3$ is the cyclic group of order twelve generated by $\tau \epsilon w$.
\end{proof}

\begin{prop}\Label{prop:I32stab}
The stabilizer in $\Gamma$ of $D'(\I^3_2)$ is 
\[\Stab{D'(\I^3_2)}=\{e,\epsilon w, w \check{\sigma}^{-1}w, w \check{\sigma}^{-1} \epsilon^3, \sigma \epsilon^3 w, \sigma\epsilon^2\}\]
and is isomorphic to $\mathfrak{S}_3$.
\end{prop}
\begin{proof}
Note that $\Stab{D'(\I^3_2)}=\Stab{\I^3_2}$.  Fix an ordering of $\I^3_2$.  Since a parabolic subgroup is its own normalizer, the group which preserves the ordering is $\Gamma\cap P_0 \cap \conj{w} \cap \conj{\sigma w}$.  The proof of Proposition~\ref{prop:I21stab} shows that $\Gamma \cap P_0 \cap \conj{w}=\langle \epsilon \rangle$ and one can easily check that $\langle \epsilon \rangle  \cap \conj{\sigma w}=\{e\}$.  Thus $\Stab{\I^3_2}$ is isomorphic to a subgroup of $\mathfrak{S}_3$.  One checks that $\{e,\epsilon w, w \check{\sigma}^{-1}w, w \check{\sigma}^{-1} \epsilon^3, \sigma \epsilon^3 w, \sigma\epsilon^2\}$ is a set of six distinct elements in $\Stab{\I^3_2}$.  Thus this set is exactly $\Stab{\I^3_2}$ and is isomorphic to $\mathfrak{S}_3$.
\end{proof}

\begin{prop}\Label{prop:I33stab}
The stabilizer in $\Gamma$ of $D'(\I^3_3)$ is trivial. 
\end{prop}
\begin{proof}
Since $|\QQ{v_0}{v_w}|^2=1$, $|\QQ{v_0}{v_{\tau \sigma w}}|^2=1$, and $|\QQ{v_w}{v_{\tau \sigma w}}|^2=2$, if $\gamma \in \Stab{\I^3_3}$, then $\lsp{\gamma}\I^2_1=\I^2_1$ or $\lsp{\gamma}\I^2_1=\{P_0,\conj{\tau \sigma w}\}$.  Thus $\Stab{\I^3_3}= \Gamma_1 \cup \Gamma_2$,  where 
\begin{align*}
\Gamma_1&=\left\{\left.\gamma \in \Stab{\I^2_1}\;\right|\;\conj{\gamma \tau \sigma w}=\conj{\tau \sigma w}\right\} \quad \text{and}\\
\Gamma_2&=\left\{\left.\gamma \in \tau \sigma\cdot\Stab{\I^2_1}\;\right|\;\conj{\gamma \tau \sigma w}=\conj{w}\right\}.
\end{align*}
One checks using Proposition~\ref{prop:I21stab} that $\Gamma_1=\{e\}$ and $\Gamma_2=\emptyset$.
\end{proof}

\begin{prop}\Label{prop:I41stab}
The stabilizer in $\Gamma$ of $D'(\I^4_1)$ is trivial.
\end{prop}
\begin{proof}
Note that $\Stab{D'(\I^4_1)}=\Stab{\I^4_1}$.  The set of isotropic vectors associated to $\I^3_1 \subset \I^4_1$ span a two dimensional subspace of $\C^3$ while the isotropic vectors associated to any other order three subset of $\I^4_1$ span all of $\C^3$.  Thus if $\gamma \in \Stab{\I^4_1}$ then $\lsp{\gamma} \I^3_1=\I^3_1$ and $\conj{\gamma \sigma w}=\conj{\sigma w}$  Therefore one calculates from Proposition~\ref{prop:I31stab} that $\Stab{\I^4_1}=\Stab{\I^3_1}\cap \conj{\sigma w}=\{e\}$.   \end{proof}

\begin{prop}\Label{prop:I42stab}
The stabilizer in $\Gamma$ of $D'(\I^4_2)$ is cyclic of order two generated by $\epsilon w$.
\end{prop}
\begin{proof}
Explicit computation show that 
\begin{align*}
{}^{\epsilon w}\I^3_3&=\{P_0,\conj{w},\conj{w^{-1}\tau\check{\sigma} w}\},\\
{}^{\sigma \tau w \tau \sigma^{-1}}\I^3_3&=\{P_0,\conj{\tau \sigma w},\conj{w^{-1}\tau\check{\sigma} w}\},\quad \text{and}\\ 
{}^{\sigma \tau w\tau \sigma^{-1} w}\I^3_3&=\{\conj{w},\conj{\tau \sigma w},\conj{w^{-1}\tau\check{\sigma} w}\}.
\end{align*}
Then $\Stab{\I^4_2}=\Gamma_1\cup\Gamma_2\cup\Gamma_3\cup\Gamma_4$, where 
\begin{align*}
\Gamma_1&=\left\{\left.\gamma \in \Stab{\I^3_3}\;\right|\;\conj{\gamma w^{-1}\tau\check{\sigma} w}=\conj{w^{-1}\tau\check{\sigma} w}\right\},\\
\Gamma_2&=\left\{\left.\gamma \in \epsilon w \cdot\Stab{\I^3_3}\;\right|\;\conj{\gamma w^{-1}\tau\check{\sigma} w}=\conj{\tau \sigma w}\right\},\\
\Gamma_3&=\left\{\left.\gamma \in \sigma \tau w \tau \sigma^{-1}\cdot\Stab{\I^3_3}\;\right|\;\conj{\gamma w^{-1}\tau\check{\sigma} w}=\conj{w}\right\}, \quad \text{and}\\
\Gamma_4&=\left\{\left.\gamma \in \sigma \tau w \tau \sigma^{-1} w\cdot\Stab{\I^3_3}\;\right|\;\conj{\gamma w^{-1}\tau\check{\sigma} w}=P_0\right\}.
\end{align*}
By Proposition~\ref{prop:I33stab}, $\Stab{\I^3_3}$ is trivial.  It is easy to check that $\Gamma_1=\{e\}$, $\Gamma_2=\{\epsilon w\}$, and $\Gamma_3=\Gamma_4=\emptyset$.    
\end{proof}
\begin{prop}\Label{prop:I5stab}
The stabilizer group of $D'(\I^5)$ is the cyclic group of order two generated by $\sigma \epsilon^2$.
\end{prop}
\begin{proof}
Note that $\Stab{D'(\I^5)}=\Stab{\I^5}$.  With the exception of $P_0$, for every $P\in \I^5$, there exists a $Q \in \I^5$ such that $|\QQ{v_P}{v_Q}|^2=2$.  Therefore if $\gamma \in \Stab{\I^5}$, then $\conj{\gamma}=P_0$.  This implies that \begin{align*}
\Stab{\I^5}&=\Gamma_{P_0} \cap \Stab{\{\conj{w},\conj{\tau w}, \conj{\sigma w},\conj{\tau \sigma w}\}}\\      
&=\Gamma_{P_0} \cap \lsp{w \tau^{-1}\epsilon w \sigma^{-1}\tau^{-1}}\Stab{\I^4_2}\\
&=\Gamma_{P_0} \cap \lsp{w \tau^{-1}\epsilon w \sigma^{-1}\tau^{-1}}\{e,\epsilon w\} &&\text{from Proposition~\ref{prop:I42stab}.}
\end{align*}
One easily checks that this intersection is $\{e,\sigma \epsilon^2\}$.
\end{proof}

\begin{prop}\Label{prop:I8stab}
The stabilizer in $\Gamma$ of $D'(\I^8)$ is the group of order $32$ generated by $\xi^2$ and $\epsilon w$ given below:
\begin{align*}
\Stab{D'(\I^8)}=\{&e,\xi^2,\xi^4,\xi^6,\epsilon w, \epsilon w \xi^2 ,\epsilon w \xi^4 ,\epsilon w \xi^6, \xi^2 \epsilon w ,\xi^2 \epsilon w \xi^2,\xi^2 \epsilon w \xi^4,\\&\xi^2 \epsilon w \xi^6,\xi^4 \epsilon w ,\xi^4 \epsilon w \xi^2,\xi^4 \epsilon w \xi^4,\xi^4 \epsilon w \xi^6,\xi^6 \epsilon w ,\xi^6 \epsilon w \xi^2,\\&\xi^6 \epsilon w \xi^4,\xi^6 \epsilon w \xi^6,\epsilon w \xi^2\epsilon w ,\epsilon w \xi^4\epsilon w ,\epsilon w \xi^6\epsilon w,\xi^2 \epsilon w \xi^2\epsilon w,\\&\xi^2 \epsilon w \xi^4\epsilon w,\xi^2 \epsilon w \xi^6\epsilon w, \xi^4 \epsilon w \xi^2\epsilon w,\xi^4 \epsilon w \xi^4\epsilon w,\xi^4 \epsilon w \xi^6\epsilon w,\\& \xi^6 \epsilon w \xi^2\epsilon w,\xi^6 \epsilon w \xi^4\epsilon w,\xi^6 \epsilon w \xi^6\epsilon w\}
\end{align*}
It is isomorphic to the group of order $32$ with Hall-Senior number $31$ \textup{\cite{Ha}}.
\end{prop}

\begin{proof}
Note that $\Stab{D'(\I^8)}=\Stab{\I^8}$.  Consider $\I^2_1 \subset \I^8$.  Since the set of isotropic vectors associated to $\I^2_1$ is 1-bounded, if $\gamma \in \Stab{\I^8}$, then the set of isotropic vectors associated to $\lsp{\gamma} \I^2_1$ is 1-bounded.  There are 16 such subsets $\I =\lsp{\gamma_\I} \I^2_1 \subset \I^8$ with this property.  Then 
\[\Stab{\I^8}=\coprod_{\substack{\I \subset \I^8\\\I \text{ is $\Gamma$-conjugate to }\I^2_1}} \Gamma_\I,\]
where 
\[\Gamma_\I=\left\{\left.\gamma \in \gamma_\I \cdot \Stab{\I^2_1}\; \right|\; \lsp{\gamma}\{\I^8 \setminus \I^2_1\}=\{\I^8 \setminus \I\}\right\}.\]
One can compute that each $\Gamma_\I$ has exactly two elements.  Explicit computation of each $\Gamma_\I$ gives the 32 elements listed.   

Fix an ordering $[\I^8]$.  This induces a homomorphism $\psi: \Stab{\I^8} \to \mathfrak{S}_8$.  The subgroup of $\Stab{\I^8}$ which fixes the ordering is $\ker(\psi) =\Gamma \cap \left(\bigcap_{P \in \I^8} P \right)$.
One calculates that this intersection is trivial.  Thus $\Stab{\I^8}$ is isomorphic to its image $\psi(\Stab{\I^8})$ in $\mathfrak{S}_8$.  In particular, non-trivial elements of  $\psi(\Stab{\I^8})$ written in disjoint cycle notation consist of 2, 4, or 8-cycles.  Thus the exponent of $\Stab{\I^8}$ is 8.  One computes that the center of $\Stab{\I^8}$ is cyclic of order four.  There are only three groups of order 32 with exponent 8 and center which is cyclic of order four, namely $\G_{26}$, $\G_{31}$, and $\G_{32}$ in the lists given in \cite{Ha}.  These are distinguished by the number of conjugacy classes of maximal elementary abelian subgroups, that is, subgroups which are isomorphic to $(\Z/2\Z)^r$ for some $r$.  $\G_{26}$, $\G_{31}$, and $\G_{32}$ have one, two, and three conjugacy classes of maximal elementary abelian subgroups, respectively.  One calculates that $\Stab{\I^8}$ has two conjugacy classes of maximal elementary abelian subgroups.  
\end{proof}
\end{section}

\begin{section}{Structure of the spine}\Label{sec:cellstructure}
\begin{subsection}{Cell structure and action of the stabilizers}\Label{sec:cell}
From the previous computations, one can get a very explicit description of the cellular structure of the spine.  For example, to find the strongly admissible sets of order three that are on the boundary of $\I^2_1$, one needs to find all rational parabolic subgroups $Q$ such that $\I=\I^2_1 \cup \{Q\}$ is $\Gamma$-conjugate to either $\I^3_1,\ \I^3_2,$ or $\I^3_3$ and hence reduces to a simple calculation of the $\mathcal{Q}$ inner-products of $v_Q$ with $(1,0,0)^t$ and $(0,0,1)^t$.

Table~\ref{tab:sasstab} summarizes the strongly admissible sets up to $\Gamma$-conjugacy and gives their stabilizers.
\begin{table}
\caption{Stabilizer groups of strongly admissible sets}
\Label{tab:sasstab}
\begin{tabular}{|c|c|c|}
\hline
S.A.S. & Stabilizer&Generators\\\hline
$\I^2_1$&$\Z/2\Z \times \Z/4\Z$&$\langle \epsilon w , \epsilon \rangle$\\
$\I^2_2$&$\Z/8\Z$&$\langle \xi \rangle$\\
$\I^3_1$&$\Z/12\Z$&$\langle \tau \epsilon w\rangle$\\
$\I^3_2$&$\mathfrak{S}_3$&$\langle \epsilon w,\sigma \epsilon^2 \rangle$\\
$\I^3_3$&trivial&$\langle e \rangle$\\
$\I^4_1$&trivial&$\langle e \rangle$\\
$\I^4_2$&$\Z/2\Z$&$\langle \epsilon w \rangle$\\
$\I^5$&$\Z/2\Z$&$\langle \sigma \epsilon^2\rangle$\\
$\I^8$&$\langle \epsilon w, \xi^2\rangle$&$\G_{31}$\\
\hline
\end{tabular}
\end{table}

\begin{figure}
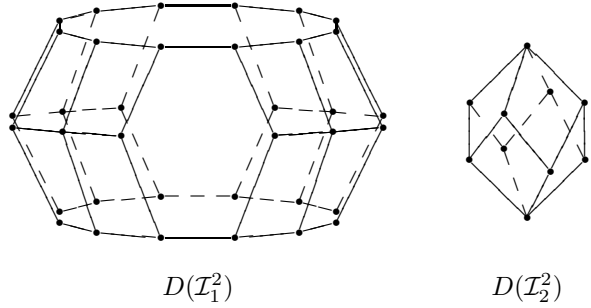

\[\begin{array}{c@{\hspace{0.4in}}c}
{\def\objectstyle{\scriptscriptstyle}
\xy /r0.5in/:="a", +(0,1)="b",+(0,1)="c","a",
{\xypolygon12"A"{~:{(1.5,0):(0,.15)::}~>{}{\bullet}}},
"b",
{\xypolygon12"B"{~:{(2,0):(0,.075)::}~>{}}},
"b",
{\xypolygon8"D"{~:{(2.1,0):(0,.075)::}~>{}{\bullet}}},
"c",
{\xypolygon12"C"{~:{(1.5,0):(0,.15)::}{\bullet}}},
"B2"*{\bullet},"B5"*{\bullet},"B8"*{\bullet},"B11"*{\bullet},
"A1";"D1"**@{--},"A2";"B2"**@{--},"A3";"D2"**@{--},"A4";"D3"**@{--},
"A5";"B5"**@{--},"A6";"D4"**@{--},"A7";"D5"**@{-},"A8";"B8"**@{-},
"A9";"D6"**@{-},"A10";"D7"**@{-},"A11";"B11"**@{-},"A12";"D8"**@{-},
"C1";"D1"**@{-},"C2";"B2"**@{--},"C3";"D2"**@{--},"C4";"D3"**@{--},
"C5";"B5"**@{--},"C6";"D4"**@{-},"C7";"D5"**@{-},"C8";"B8"**@{-},
"C9";"D6"**@{-},"C10";"D7"**@{-},"C11";"B11"**@{-},"C12";"D8"**@{-},
"D5";"B8"**@{-},"D6";"B8"**@{-},"D7";"B11"**@{-},"B11";"D8"**@{-},
"D5";"B8"**@{-},"D6";"B8"**@{-},"D7";"B11"**@{-},"B11";"D8"**@{-},
"D4";"B5"**@{--},"B5";"D3"**@{--},"D2";"B2"**@{--},"B2";"D1"**@{--},
"A7";"A6"**@{--},"A6";"A5"**@{--},"A5";"A4"**@{--},"A3";"A4"**@{--},"A3";"A2"**@{--},"A2";"A1"**@{--},"A1";"A12"**@{--},
"A7";"A8"**@{-},"A9";"A8"**@{-},"A9";"A10"**@{-},"A10";"A11"**@{-},"A11";"A12"**@{-}
\endxy}&
{
\def\objectstyle{\scriptscriptstyle}
\xy /r0.3in/:="a", +(0,1)="b",+(0,1)="c",+(0,1)="d",
"a"*{\bullet},"d"*{\bullet},
"b",{\xypolygon4"B"{~={45}~:{(1,-0.14):(0.4,.2)::}~>{}{\bullet}}},
"c",{\xypolygon4"C"{~={0}~:{(1,0):(.4,.2)::}~>{}{\bullet}}},
"C1";"B1"**@{-},"C1";"B4"**@{-},
"C2";"B2"**@{--},"C2";"B1"**@{--},
"C3";"B3"**@{-},"C3";"B2"**@{--},
"C4";"B4"**@{-},"C4";"B3"**@{-},
"B1";"a"**@{-},"B2";"a"**@{--},"B3";"a"**@{-},"B4";"a"**@{-},
"C1";"d"**@{-},"C2";"d"**@{--},"C3";"d"**@{-},"C4";"d"**@{-}
\endxy}\\
\\
D(\I^2_1)&D(\I^2_2)
\end{array}\]
\caption{The sets associated to the two $\Gamma$-conjugacy classes of strongly admissible sets of order two are shown here.  $D(\I^2_1)$ is homeomorphic to a polytope with dodecagon, hexagon and quadrilateral faces, while $D(\I^2_2)$ has only quadrilateral faces.}
\Label{fig:round}
\end{figure}
\begin{table}

\caption{Incidence types}
\label{tab:incidence}
\begin{tabular}{|c|cc|ccc|cc|cc|}
\hline
 & $\I^2_1$&$\I^2_2$&$\I^3_1$&$\I^3_2$&$\I^3_3$&$\I^4_1$&$\I^4_2$&$\I^5$&$\I^8$\\\hline
$\I^2_1$&$\ast$&$\ast$&3&3&2&5&4&8&16\\
$\I^2_2$&$\ast$&$\ast$&0&0&1&1&2&2&8\\\hline
$\I^3_1$&2&0&$\ast$&$\ast$&$\ast$&1&0&2&0\\
$\I^3_2$&4&0&$\ast$&$\ast$&$\ast$&1&0&2&0\\
$\I^3_3$&12&8&$\ast$&$\ast$&$\ast$&2&4&6&32\\\hline
$\I^4_1$&40&8&12&6&2&$\ast$&$\ast$&4&0\\
$\I^4_2$&16&8&0&0&2&$\ast$&$\ast$&1&16\\\hline
$\I^5$&32&8&12&6&3&2&1&$\ast$&$\ast$\\
$\I^8$&4&2&0&0&1&0&1&$\ast$&$\ast$\\
\hline
\end{tabular}
\end{table}
The incidence table is given in Table~\ref{tab:incidence}, where the entry below the diagonal means that each column cell has that many row cells in its boundary, and the entry above the diagonal means the column cell appears in the boundary of this many row cells.  The entries below the diagonal can be read off from Figures~\ref{fig:order4}-\ref{fig:order2}, while the entries above the diagonal can be easily computed from Proposition~\ref{prop:Qmatrix}, since the $\Gamma$-conjugacy class of a strongly admissible can be distinguished by the pairwise $\mathcal{Q}$-inner products of its associated isotropic vectors, except to distinguish $\I^3_1$ and $\I^3_2$, we must also compute the dimension of the span of their isotropic vectors.  

Define a map ${\mathfrak Q}:\{\text{ordered subsets of $\Par$}\} \to \text{Mat}(\Z)$ as follows:  Given an ordered subset $[\I]=(P_1, \ldots, P_n) \subset \Par$, let $\QQQ{\I}$ be the $n \times n$ matrix whose $(i,j)$-component is $\left|\QQ{v_{P_i}}{v_{P_j}}\right|^2$.

\begin{prop}\label{prop:Qmatrix}
Let $\I$ be strongly admissible set.  Then 
\[\QQQ{\I}=\begin{cases} \mat{0&1\\1&0} & \text{if $\I$ is $\Gamma$-conjugate to $\I^2_1$,}\\
\mat{0&2\\2&0} & \text{if $\I$ is $\Gamma$-conjugate to $\I^2_2$,}\\
\mat{0&1&1\\1&0&1\\1&1&0} & \text{if $\I$ is $\Gamma$-conjugate to $\I^3_1$ or $\I^3_2$,}\\
\mat{0&1&2\\1&0&1\\2&1&0} & \text{if $\I$ is $\Gamma$-conjugate to $\I^3_3$,}\\
\mat{0&1&1&2\\1&0&1&1\\1&1&0&1\\2&1&1&0} & \text{if $\I$ is $\Gamma$-conjugate to $\I^4_1$,}\\
\mat{0&1&1&2\\1&0&2&1\\1&2&0&1\\2&1&1&0} & \text{if $\I$ is $\Gamma$-conjugate to $\I^4_2$,}\\
\mat{0&1&1&1&2\\1&0&1&2&1\\1&1&0&1&1\\1&2&1&0&1\\2&1&1&1&0}& \text{if $\I$ is $\Gamma$-conjugate to $\I^5$,} \\
\mat{0&1&2&1&1&2&1&4\\1&0&1&2&2&1&4&1\\2&1&0&1&1&4&1&2\\1&2&1&0&4&1&2&1\\1&2&1&4&0&1&2&1\\2&1&4&1&1&0&1&2\\1&4&1&2&2&1&0&1\\4&1&2&1&1&2&1&0}&\text{if $\I$ is $\Gamma$-conjugate to $\I^8$.}\end{cases}\]
\end{prop}

In Figures~\ref{fig:order4}-\ref{fig:order2}, a 0-cell that is $\Gamma$-conjugate to $\I^5$ will be denoted by $\bullet$, and a 0-cell that is $\Gamma$-conjugate to $\I^8$ will be denoted by $\circledast$.   A 1-cell that is $\Gamma$-conjugate to $\I^4_1$ will be denoted with a solid line, and a 1-cell that is $\Gamma$-conjugate to  $\I^4_2$ will be denoted with a dotted line.

The 1-cell $D(\I^4_1)$ is shown in Figure~\ref{fig:order4}.  The boundary consists of two 0-cells that are both conjugate to $D(\I^5)$.  By Proposition~\ref{prop:I41stab}, the stabilizer of $D(\I^4_1)$ is trivial.

\begin{figure}
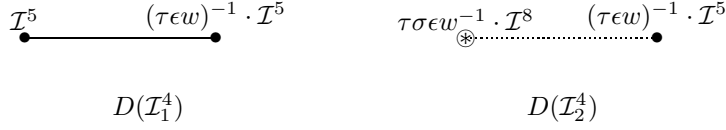

\[\begin{array}{c @{\hspace{0.5in}}c}
{\xy /l.5in/:,(-1,0)="A","A"*{\bullet},"A"*+!D{(\tau \epsilon w)^{-1}\cdot\I^5},
(1,0)="B","B"*{\bullet},"B"*+!D{\I^5},
"A";"B"**@{-} \endxy}&{\xy /l.5in/:,(-1,0)="A","A"*{\bullet},"A"*+!D{(\tau \epsilon w)^{-1}\cdot\I^5},
(1,0)="B","B"*{\circledast},"B"*+!D{\tau \sigma\epsilon w^{-1}\cdot \I^8},
"A";"B"**@{.} \endxy}\\[0.2in]
D(\I^4_1)&D(\I^4_2)
\end{array}\]
\caption{Strongly admissible sets of order 4}
\Label{fig:order4}
\end{figure}

The 1-cell $D(\I^4_2)$ is shown in Figure~\ref{fig:order4}.  The boundary consists of two 0-cells (one that is conjugate to $D(\I^5)$ and one that is conjugate to $D(\I^8)$).  By Proposition~\ref{prop:I42stab}, the stabilizer of $D(\I^4_2)$ is isomorphic to $\Z/2\Z$ and is generated by $\epsilon w$.  This acts on the cell by fixing it pointwise.

The 2-cell $D(\I^3_1)$ is shown in Figure~\ref{fig:order3}.  Its boundary consists of twelve 1-cells (all of which are conjugate to $D(\I^4_1)$) and twelve 0-cells (all of which are conjugate to $D(\I^5)$).  By Proposition~\ref{prop:I31stab} the stabilizer of $D(\I^3_1)$ is isomorphic to $\Z/12\Z$ and is generated by $\tau \epsilon w$.  It is easily checked that $\tau \epsilon w$ acts on the figure by rotation by $\frac{\pi}{6}$ about the first contact $z(\I^3_1)$ for $D(\I^3_1)$.  

\begin{figure}
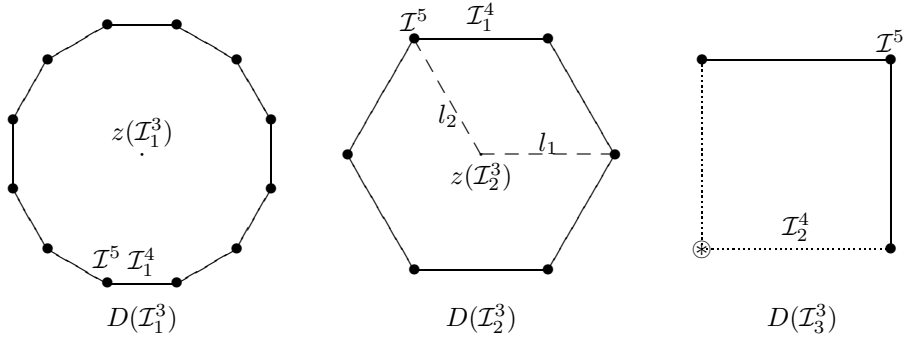

\[\begin{array}{c@{\hspace{0.4in}}c@{\hspace{0.4in}}c}
{\xy /l.7in/:,{\xypolygon12"A"{\bullet}}, "A0"*{\cdot}, "A3"*++!D{\I^5}, "A0"+R*+!D{z(\I^3_1)}, "A3";"A4" **@{} ? *+!D{\I^4_1}
\endxy}&
{\xy /l.7in/:,{\xypolygon6"B"{\bullet}}, "B0"*{\cdot}, "B0"*+!U{z(\I^3_2)}, "B0"; "B4"**@{--} ? *!D{l_1}, "B0";"B6" **@{--} ? *+!U{l_2}, "B6"*+!D{\I^5}, "B6";"B5" **@{} ? *+!D{\I^4_1}
\endxy}&
{\xy /l.7in/:,{\xypolygon4"C"{~>{}}}, "C3"*{\bullet}, "C3"*+!D{\I^5},"C2"*{\bullet},"C4"*{\bullet},"C1"*{\circledast}, "C1";"C4"**@{.}, "C1";"C2"**@{.} ?*+!D{\I^4_2}, "C2";"C3" **@{-}, "C3";"C4" **@{-}
\endxy}\\[0.7in]
D(\I^3_1)& D(\I^3_2)&D(\I^3_3)
\end{array}\]
\caption{Strongly admissible sets of order 3}
\Label{fig:order3}
\end{figure}

The 2-cell $D(\I^3_2)$ is shown in Figure~\ref{fig:order3}.  The boundary consists of six 1-cells (all of which are conjugate to $D(\I^4_1)$) and six 0-cells (all of which are conjugate to $D(\I^5)$).  The point $z(\I^3_2)$ is the first contact for $\I^3_2$.  The lines $l_j$ represent the gradient flows of the $ f_0$ function restricted to $D(\I^3_2)$ from the first contact to two of the  vertices of $D(\I^3_2)$.  By Proposition~\ref{prop:I32stab}, the stabilizer of $D(\I^3_2)$ is isomorphic to $\mathfrak{S}_3$ generated by $\epsilon w$ and $\sigma^2 \epsilon^2$.  One can check that $\epsilon w$ acts as reflection about $l_1$ and $\sigma \epsilon^2$ acts as reflection about $l_2$.
 
The 2-cell $D(\I^3_3)$ is shown in Figure~\ref{fig:order3}.  The boundary consists of four 1-cells (two that are conjugate to $D(\I^4_1)$ and two that are conjugate to $D(\I^4_2)$) and four 0-cells (three that are conjugate to $D(\I^5)$ and one that is conjugate to $D(\I^8)$).  By Proposition~\ref{prop:I33stab}, the stabilizer of $D(\I^3_3)$ is trivial.

See Figure~\ref{fig:order2} for a picture of the face relations for $D(\I^2_1)$.  The outside of the figure is the dodecagon bottom face seen in Figure~\ref{fig:round}.  The boundary of $D(\I^2_1)$ consists of twenty-two 2-cells (two that are conjugate to $D(\I^3_1)$, four that are conjugate to $D(\I^3_2)$, and sixteen that are conjugate to $D(\I^3_3)$), fifty-six 1-cells (forty that are conjugate to $D(\I^4_1)$ and sixteen that are conjugate to $D(\I^4_2)$), and thirty-six 0-cells (thirty-two that are conjugate to $D(\I^5)$ and four that are conjugate to $D(\I^8)$).  By Proposition~\ref{prop:I21stab} the stabilizer of $D(\I^2_1)$ is isomorphic to $\Z/4\Z \times \Z/2\Z$ generated by $\epsilon$ and $\epsilon w$.  It is easily calculated that $\epsilon$ acts by rotating the figure by $-\frac{\pi}{2}$ and $\epsilon w$ acts by an inversion, sending the interior 12-gon to the exterior one.   
\begin{figure}
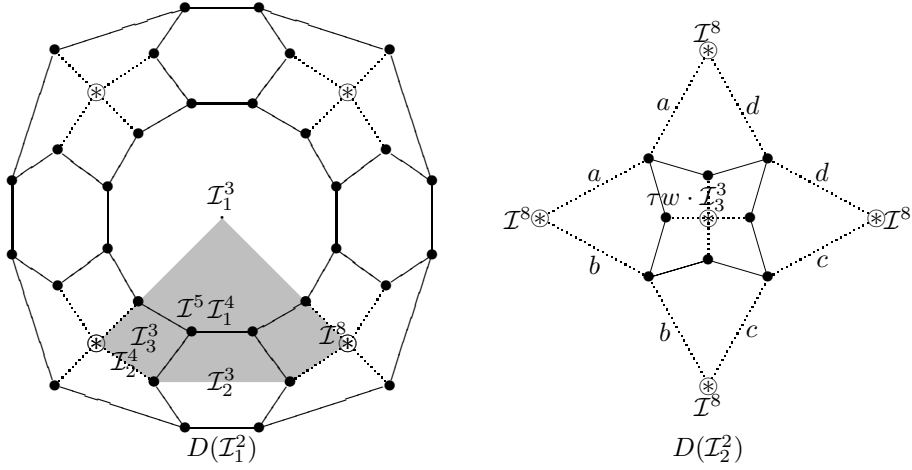

\[\begin{array}{c @{\hspace{0.3in}} c}
{\xy /l.62in/:,{\xypolygon12"A"{\bullet}}, 
{\xypolygon4"B"{~:{(1.5,0):}~>{}{\circledast}}},
{\xypolygon4"C"{~:{(2,0):}~>{}{\bullet}}},
{\xypolygon8"D"{~:{(1.5,0):}~>{}{\bullet}}},
{\xypolygon36"E"{~:{(1.8,0):}~>{}}},
@={"A0","B1","D2","D3","B2","A0"},{0*[lightgrey]\xypolyline{*}},
"A0"*{\cdot},
"A2"*{\bullet},"A3"*{\bullet},"A4"*{\bullet},"A5"*{\bullet},
"A2";"A3"**@{-},"A3";"A4"**@{-},"A4";"A5"**@{-},
"B1"*{\circledast},"B2"*{\circledast},
"D2"*{\bullet},"D3"*{\bullet},
"A3"*++!D{\I^5}, "A0"+R*+!D{\I^3_1}, "A3";"A4" **@{} ? *+!D{\I^4_1}, 
"B1";"A2"**@{.},"B2";"A5"**@{.}, "B3";"A8"**@{.},"B4";"A11"**@{.}, 
"B1";"C1"**@{.}, "B2";"C2"**@{.}, "B3";"C3"**@{.}, "B4";"C4"**@{.},
"A3";"D2"**@{-},"A4";"D3"**@{-},"A6";"D4"**@{-},"A7";"D5"**@{-},"A9";"D6"**@{-},
"A10";"D7"**@{-},"A12";"D8"**@{-},"A1";"D1"**@{-},
"B1";"D2"**@{.}?*{\I^4_2},"B2";"D4"**@{.},"B3";"D6"**@{.},"B4";"D8"**@{.},
"B1";"D1"**@{.},"B2";"D3"**@{.},"B3";"D5"**@{.},"B4";"D7"**@{.},
"C1";"E2"**@{-},"C2";"E18"**@{-},"C3";"E27"**@{-},"C4";"E36"**@{-},
"C1";"E9"**@{-},"C2";"E11"**@{-},"C3";"E20"**@{-},"C4";"E29"**@{-},
"E11";"E9"**@{-},"E18";"E20"**@{-},"E29";"E27"**@{-},"E36";"E2"**@{-},
"D1";"E2"**@{-},"D4";"E18"**@{-},"D6";"E27"**@{-},"D8";"E36"**@{-},
"D2";"E9"**@{-},"D3";"E11"**@{-},"D5";"E20"**@{-},"D7";"E29"**@{-},
"E2"*{\bullet},"E18"*{\bullet},"E27"*{\bullet},"E36"*{\bullet},
"E9"*{\bullet},"E11"*{\bullet},"E20"*{\bullet},"E29"*{\bullet},
"D3";"D2"**@{}?*{\I^3_2},"A3";"B1"**@{}?!{"A2";"D2"}*{\I^3_3},
"B2"*!RD{\I^8},
\endxy}&

{\xy /l.22in/:,{\xypolygon4"A"{~={90}~>{}{\bullet}}},
{\xypolygon4"B"{~:{(2,0):}~>{}{\bullet}}},
{\xypolygon4"C"{~:{(4,0):}~={90}~>{}{\circledast}}},
"A0"*{\circledast},
@={"A1","B1","A4","B4","A3","B3","A2","B2","A1","B1"},
{\xypolyline{}},
"A1";"A3"**@{.},
"A2";"A4"**@{.},
"A3";"A4"**@{}?*{\tau w \cdot \I^3_3},
"C1"*+!U{\I^8},"C2"*+!L{\I^8},"C3"*+!D{\I^8},"C4"*+!R{\I^8},
"C1";"B1"**@{.}?*+!R{b},
"B1";"C4"**@{.}?*+!U{b},
"C4";"B4"**@{.}?*+!D{a},
"B4";"C3"**@{.}?*+!R{a},
"C3";"B3"**@{.}?*+!L{d},
"B3";"C2"**@{.}?*+!D{d},
"C2";"B2"**@{.}?*+!U{c},
"B2";"C1"**@{.}?*+!L{c}
\endxy}\\[.5in]
D(\I^2_1)& D(\I^2_2)
\end{array}\]
\caption{Strongly admissible sets of order 2}
\Label{fig:order2}
\end{figure}

The other 3-cell $D(\I^2_2)$ is shown in Figure~\ref{fig:order2}.  The labels $\{a,b,c,d\}$ are added to show the identifications that need to be made.  The boundary of $D(\I^2_2)$ consists of eight 2-cells (all of which are conjugate to $D(\I^3_3)$), sixteen 1-cells (eight that are conjugate to $D(\I^4_1)$ and  eight that are conjugate to $D(\I^4_2)$) and ten 0-cells (eight that are conjugate to $D(\I^5)$ and two that are conjugate to $D(\I^8)$).  By Proposition~\ref{prop:I22stab}, the stabilizer of $D(\I^2_2)$ is isomorphic to $\Z/8\Z$ generated by $\xi$.  It is easily checked that $\xi$ acts on the figure by the composition of an inversion sending the exterior point to the interior point and a rotation of $\frac{\pi}{4}$.
\end{subsection}

\begin{subsection}{The subdivision}
From the explicit description of the cells and $\Gamma$-action, it is straightforward to subdivide the spine so that the stabilizer of each cell fixes the cell pointwise. 

First consider $D(\I^2_1)$.  It must be divided into eight 3-cells that are conjugate via $\Stab{\I^2_1}$.  One of the 3-cells $X$ is represented by the shaded region in Figure~\ref{fig:order2}.  It is easy to see that the boundary of $X$ consists of two 2-cells, $C_1$ and $C_2$, that are $\Gamma$-conjugate to $D(\I^3_3)$, half of the 2-cell $D(\I^3_2)$ $B$, one-fourth of the 2-cell $D(\I^3_1)$, $A$, and three new faces, $E_1$, $E_2$, and $D$ that lie inside $D(\I^2_1)$ and such that $E_1$ and $E_2$ are $\Gamma$-conjugate.  Next we turn to subdividing the boundary faces of $X$.  From the description of the action of the $\Stab{\I^3_1}$ on $D(\I^3_1)$, it follows that $A$ must be subdivided into three $\Gamma$-conjugate 2-cells $A_1$, $A_2$, and $A_3$.  Similarly, because of the action of $\Stab{\I^3_2}$ on $D(\I^3_2)$, $B$ must be subdivided into three $\Gamma$-conjugate 2-cells $B_1$, $B_2$, and $B_3$.  One can compute that the stabilizers of $E_i$ are trivial, and that the stabilizer of $D$ fixes $D$ pointwise.  The 1-cells of $X$, both the ones that are $\Gamma$-conjugate to $D(\I^4_j)$ $(j=1,2)$ and the new ones that are introduced for the subdivision, do not need to be subdivided because the stabilizers are either trivial, or the stabilizer acts on the 1-cell by fixing it pointwise.  This yields $X$ in Figure~\ref{fig:subdivided}.

Similarly, $D(\I^2_2)$ must be subdivided into eight 3-cells that are conjugate via $\Stab{\I^2_2}$.  One of the 3-cells $Y$ can be viewed in Figure~\ref{fig:order2} by taking a face $C_1$ that is $\Gamma$-conjugate to $D(\I^3_3)$ and looking at the set of points of $D(\I^2_2)$ that would flow to $C_1$ under the gradient flow for $-f_{\I^2_2}$ and adding the first contact point of $\I^2_2$.  Then the boundary of $Y$ consists of five $2$-cells, $C_1$, $F_1$, $F_2$, $G_1$, and $G_2$ such that $F_1$ is $\Gamma$-conjugate to $F_2$ and $G_1$ is $\Gamma$-conjugate to $G_2$.  The calculations above show that the stabilizer of $C_1$ is trivial, and it is easy to check that the stabilizers of $F_1$, $F_2$, $G_1$, and $G_2$ are trivial as well.  Hence the 2-cells do not need to be subdivided.  Similarly, one checks that the stabilizers of the 1-cells are either trivial or fix the 1-cell pointwise.  This yields $Y$ in Figure~\ref{fig:subdivided}.
\begin{figure}
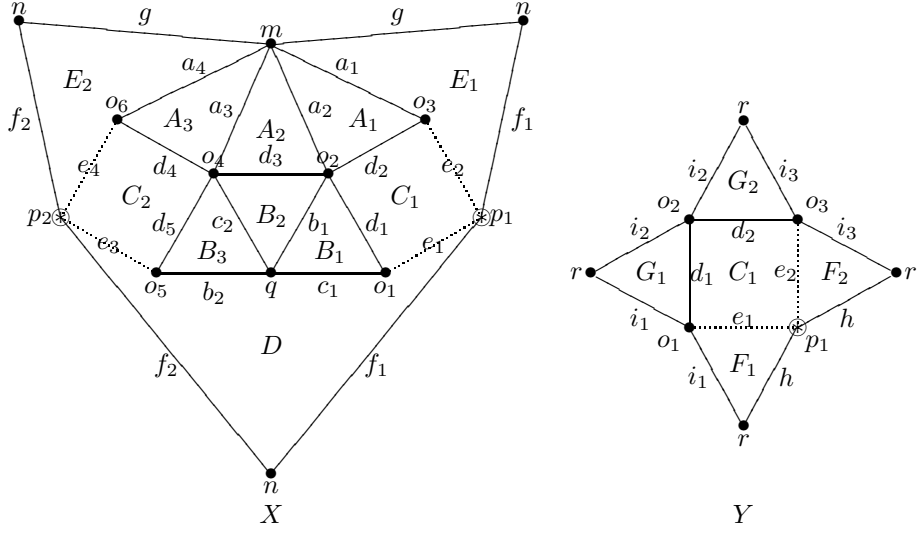

\[\begin{array}{c @{\hspace{0.1in}} c}
{\xy /l.6in/:="c",{\xypolygon6"A"{~>{}}}, 
{\xypolygon12"B"{~:{(1.9,0):}~>{}}},
"c",+(2.2,-2.2)="D","D"*{\bullet}*+!D{n},
"c",+(-2.2,-2.2)="E","E"*{\bullet}*+!D{n},
"c",+(0,-2)="F",
"c",+(0,1.75)="G",
"F";"G"**@{}?(0.2)*{A_2}?(0.4)*{B_2}?(0.7)*{D},
"D";"B12"**@{-}?*+!R{f_2},"E";"B7"**@{-}?*+!L{f_1},
"D";"F"**@{-}?*+!D{g},"E";"F"**@{-}?*+!D{g},
"D";"B11"**@{}?(0.6)*{E_2},"E";"B8"**@{}?(0.6)*{E_1},
"A1";"B11"**@{}?*{C_2},"A4";"B8"**@{}?*{C_1},
"B8";"B11"**@{}?(0.8)*{A_3}?(0.2)*{A_1},
"A1"*{\bullet}*+!U{o_5},
"A4"*{\bullet}*+!U{o_1},
"A5"*{\bullet}*+!D{o_2},
"A6"*{\bullet}*+!D{o_4},
"A0"*{\bullet},"A0"*+!U{q},
"A1";"A0"**@{-}?*+!U{b_2},
"A4";"A0"**@{-}?*+!U{c_1},
"A5";"A0"**@{-}?*+!L{b_1},
"A6";"A0"**@{-}?*+!R{c_2},
"A6";"A1"**@{-}?*+!R{d_5},
"A5";"A6"**@{-}?*+!D{d_3},
"A4";"A5"**@{-}?*+!L{d_1},
"B12"*{\circledast}*+!R{p_2},
"B11"*{\bullet}*+!D{o_6},
"B7"*{\circledast}*+!L{p_1},
"B8"*{\bullet}*+!D{o_3},
"B7";"B12"**@{}?(0.36)*+++!U{B_1}?(0.64)*+++!U{B_3},
"B12";"B11"**@{.}?*{e_4},
"B7";"B8"**@{.}?*{e_2},
"B12";"A1"**@{.}?*{e_3},
"B11";"A6"**@{-}?*+!U{d_4},
"B7";"A4"**@{.}?*{e_1},
"B8";"A5"**@{-}?*+!U{d_2},
"F"*{\bullet}*+!D{m},
"G"*{\bullet}*+!U{n},
"F";"A6"**@{-}?*+!R{a_3},
"F";"A5"**@{-}?*+!L{a_2},
"F";"B11"**@{-}?*+!D{a_4},
"F";"B8"**@{-}?*+!D{a_1},
"G";"B12"**@{-}?*+!U{f_2},
"G";"B7"**@{-}?*+!U{f_1}
\endxy}&
{\xy /l.4in/:,{\xypolygon4"B"{~>{}}},
{\xypolygon4"C"{~:{(2,0):}~={90}~>{}{\bullet}}},
"B0"*{C_1},
"B1"*{\bullet},
"B1"*+!RU{o_1},
"B2"*{\circledast},
"B2"*+!LU{p_1},
"B3"*{\bullet},
"B3"*+!LD{o_3},
"B4"*{\bullet},
"B4"*+!RD{o_2},
"B1";"B2"**@{.}?*!D{e_1},
"B3";"B2"**@{.}?*!R{e_2},
"B3";"B4"**@{-}?*!U{d_2},
"B1";"B4"**@{-}?*!L{d_1},
"C1"*+!U{r},"C2"*+!L{r},"C3"*+!D{r},"C4"*+!R{r},
"C1";"B1"**@{-}?*+!R{i_1},
"B1";"C4"**@{-}?*+!U{i_1},
"C4";"B4"**@{-}?*+!D{i_2},
"B4";"C3"**@{-}?*+!R{i_2},
"C3";"B3"**@{-}?*+!L{i_3},
"B3";"C2"**@{-}?*+!D{i_3},
"C2";"B2"**@{-}?*+!U{h},
"B2";"C1"**@{-}?*+!L{h},
"C3";"C1"**@{}?(0.8)*{F_1},
"C3";"C1"**@{}?(0.2)*{G_2},
"C2";"C4"**@{}?(0.8)*{G_1},
"C2";"C4"**@{}?(0.2)*{F_2}
\endxy}\\[0.7in]
X& Y
\end{array}\]
\caption{Subdivided cells}
\Label{fig:subdivided}
\end{figure}

In the figures, the labels that only differ by a subscript are conjugate under $\Gamma$.  For each $\Gamma$-conjugacy class, we fix a representative and compute the stabilizers.  The results are given in Table~\ref{tab:stab}.
\begin{table}
\caption{Representative Cells and Their Stabilizers}
\Label{tab:stab}
\begin{center}
\begin{tabular}{|c|c|c|c|}\hline
Cell & Dimension & Stabilizer & Generators \\\hline
$X$ & 3& trivial & $\langle e \rangle$ \\ 
$Y$ & 3&trivial & $\langle e \rangle$ \\
$A_1$ &2& trivial &$\langle e \rangle$ \\
$B_1$ &2& trivial &$\langle e \rangle$ \\
$C_1$ &2& trivial &$\langle e \rangle$ \\
$D$ &2& $\Z/2\Z$ & $\langle \epsilon w \rangle$\\
$E_1$ &2& trivial &$\langle e \rangle$ \\
$F_1$ &2& trivial &$\langle e \rangle$\\
$G_1$ &2& trivial &$\langle e \rangle$ \\
$a_1$ &1& trivial &$\langle e \rangle$ \\
$b_1$ &1& $\Z/2\Z$ & $\langle \sigma \epsilon w \sigma^{-1} \rangle$\\
$c_1$&1&$\Z/2\Z$ & $\langle\epsilon w\rangle$\\
$d_1$&1& trivial &$\langle e \rangle$ \\
$e_1$&1&$\Z/2\Z$ &  $\langle \epsilon w \rangle$\\
$f_1$&1&$\Z/2\Z$ &  $\langle\epsilon w \rangle$\\
$g$&1&$\Z/4\Z$ & $\langle \epsilon \rangle$\\
$h$&1&$\Z/4\Z$ & $\langle \xi^2 \rangle$\\
$i_1$&1&trivial &$\langle e \rangle$ \\
$m$&0&$\Z/12\Z$ & $\langle \tau \epsilon w \rangle$\\
$n$&0&$ \Z/4\Z \times \Z/2\Z$ & $\langle \epsilon w , w \rangle$\\
$o_1$&0&$\Z/2\Z$ & $\langle \epsilon w \rangle $\\
$p_1$&0&$\G_{31}$ & $\langle \epsilon w , \xi^2 \rangle$\\
$q$&0&$\mathfrak{S}_3$ &$\langle \epsilon w,\sigma \epsilon^2 \rangle$\\
$r$&0&$\Z/8$ &$\langle \xi \rangle$\\\hline
\end{tabular}
\end{center}
\end{table}
\end{subsection}
\end{section}

\begin{section}{Cohomology of $\Gamma \backslash D_0$ with local coefficients}\Label{sec:cohomology}
In this section only, we generalize to the case where $D=G/K$ is a non-compact symmetric space, where $G$ is the group of real points of a semisimple linear algebraic group $G$ defined over $\Q$.  Let $\Gamma$ be an arithmetic subgroup of $G(\Q)$.  Let $D_0\subset D$ denote a spine with CW-structure such that the $\Gamma$-stabilizer of a cell fixes the cell pointwise.

In order to set notation, we recall the definition of orbifold and the sheaf associated to the local system.  We then prove that the cohomology of $\Gamma\backslash D$ with local coefficients is isomorphic the the cohomology of $\Gamma \backslash D_0$ with local coefficients.
\begin{subsection}{Orbifolds}
The notion of an orbifold was first introduced by Satake in \cite{Sat}.  He called them V-manifolds.  Let $M$ be a Hausdorff topological space.  A \emph{local uniformizing system \textup (l.u.s.\textup)} $\{U,\Gamma_U,\phi\}$ for an open set $V\subset M$ is a collection of the following objects:
\begin{enumerate}
\item $U$: a connected open subset of $\R^n$.
\item $G$: a finite group of linear transformations of $U$ to itself such that the set of fixed points of $G$ is either all of $U$ or at least codimension 2.
\item $\phi$: a continuous $\Gamma_U$-invariant map $U \to V$ such that the induced map $\Gamma_U \backslash U \to V$ is a homeomorphism.   
\end{enumerate}
Let $V \subset V' \subset M$ be two open sets and let $\{U,\Gamma_U,\phi\}$ and $\{U',\Gamma_{U'},\phi'\}$ be local uniformizing systems for $V$ and $V'$ respectively.   An \emph{injection} $\lambda:\{U,\Gamma_U,\phi\}\to \{U',\Gamma_{U'},\phi'\}$ is a smooth injection $\lambda:U \to U'$ such that for any $\gamma \in \Gamma_U$, there exists a $\gamma' \in \Gamma_{U'}$ such that $\lambda \circ \gamma = \gamma' \circ \lambda$ and $\phi = \phi' \circ \lambda$.   

Let $\mathfrak{L}$ be a family of local uniformizing systems for open sets in $M$.  Then an open set $V\subset M$ is said to be \emph{$\mathfrak{L}$-uniformized} if there exists a local uniformizing system for $V$ in $\mathfrak{L}$.

\begin{defn}
An \emph{orbifold} is a Hausdorff space $M$ and a family $\mathfrak{L}$ of local uniformizing systems for open sets in $M$ satisfying the following conditions. 
\begin{enumerate}
\item Let $V \subset V' \subset M$ be two open sets and let $\{U,\Gamma_U,\phi\}$ and $\{U',\Gamma_{U'},\phi'\}$ be local uniformizing systems for $V$ and $V'$ respectively.  Then there exists an injection $\{U,\Gamma_U,\phi\} \to \{U',\Gamma_{U'},\phi'\}$.
\item The $\mathfrak{L}$-uniformized open sets form a basis of open sets for $M$.
\end{enumerate} 
\end{defn}

Two families of local uniformizing systems $\mathfrak{L}_1$ and $\mathfrak{L}_2$ defining an orbifold $M$ are said to be \emph{equivalent} if $\mathfrak{L}_1 \cup \mathfrak{L}_2$ satisfies the conditions above.  Equivalent families define the same orbifold structure on $M$.  Thus when talking about  family of local uniformizing systems for an orbifold $M$, we will mean a maximal family.

A \emph{smooth function} on $M$ is given locally on $\mathfrak{L}$-uniformized sets $V$ by $\Gamma_U$-invariant smooth functions on $U$.  Similarly, a \emph{smooth $p$-form}  on $M$ is given locally on $V$ by $\Gamma_U$-invariant smooth forms on $U$.

Note that a smooth manifold is an example of an orbifold where every group $\Gamma_U$ for $\{U,\Gamma_U,\phi\} \in \mathfrak{L}$ is the trivial group.  If $\tilde{M}$ is a smooth manifold and $\Gamma$ is a properly discontinuous group of automorphisms of $\tilde{M}$, then the quotient $\Gamma \backslash \tilde{M}$ has a canonical structure of an orbifold.
\end{subsection}

\begin{subsection}{$\Gamma \backslash D$ as an orbifold}
Since $\Gamma$ acts properly discontinuously on $D$, $X=\Gamma \backslash D$ has a canonical structure of an orbifold.  Let $\pi$ denote the projection $D \to X$.  

Given a connected open set $V\subset X$ with local uniformizing system $L$, we can and will identify $L$ with a triple $\{U,\Stab{U},\pi|_U\}$ where $U$ is a connected component of $\tilde{U}=\pi^{-1}(V)$.  A smooth function on $V$ is given by a $\Stab{U}$-invariant smooth function on $U$.  In particular, for an open subset $O\subset X$, a smooth function on $O$ is given by a $\Gamma$-invariant function on $\pi^{-1}(O)$. 

Since $D_0$ is a spine, there exists a $\Gamma$-equivariant deformation retract $\tilde{r}:D \to D_0$.  Then there is an induced deformation retract $r:X \to X_0 =\Gamma \backslash D_0$.
\end{subsection}  
\begin{subsection}{Locally constant orbifold (fiber) bundles over $\Gamma \backslash D$}
Let $E$ be a $\Gamma$-module with $\Gamma$-action given by $\rho:\Gamma \to \GL(E)$.  Let $D\times_\Gamma E$ denote the quotient space $D\times E/\{(x,v)\sim(\gamma\cdot x, \rho(\gamma)v)\,|\, \gamma \in \Gamma\}$.  For a $\Gamma$-invariant deformation retract $D_0 \subseteq D$, we can similarly define $D_0 \times_\Gamma E$.  This yields the following maps.
\begin{displaymath}
\xy*!<0pc,3.5pc>\xybox{\xymatrix{&D_0 \times E \ar[dl]_-{\tilde \pi_0} \ar[d] \ar[r] &D\times E\ar[d] \ar[dr]^-{\tilde \pi}& \\
D_0 \times_\Gamma E \ar[d] \ar@/^7pc/[rrr]& D_0 \ar[dl]_-{\pi_0} \ar@/_/[r]&D\ar@/_/[l]_-{\tilde r} \ar[dr]^-\pi&D \times_{\Gamma} E \ar[d]\\
{\makebox[0.4in][r]{$X_0=\Gamma \backslash D_0$}} \ar@/_/[rrr]& & &{\makebox[0.3in][l]{$\Gamma \backslash D=X$}}\ar@/_/[lll]_-{r}}}\endxy
\end{displaymath}

Suppose that $E$ is an $N$-dimensional vector space.  If $\Gamma$ is torsion-free, then $X$ is a smooth manifold and $D\times_\Gamma E$ is a flat rank-$N$ vector bundle over $X$.  If $\Gamma$ has torsion, then $X$ is an orbifold and $D\times_\Gamma E$ is called a \emph{flat orbifold bundle} over $X$.

Consider the presheaf $\mathcal{E}$ on $X$ defined as follows.  For every open set $U \subset X$, 
\begin{align*}
\mathcal{E}(U)&=\left\{f:\pi^{-1}(U) \to E\;\left | \;\substack{\text{$f$ is locally constant and}\\ f(\gamma \cdot x)=\rho(\gamma)f(x) \quad \forall \gamma \in \Gamma,\  x\in \pi^{-1}(U)}\right.\right\}
\end{align*}
with the obvious restriction maps.  Let $\E$ denote the sheafification of $\mathcal{E}$.  Similarly define the presheaf $\mathcal{E}_0$ on $X_0$ and let $\E_0$ denote the sheafification of $\mathcal{E}_0$.  Note that if  $\Gamma$ is torsion-free, $\E$ is the sheaf associated to the local system $X$ defined by $(E,\rho)$.  We will extend this terminology to $X$, respectively $X_0$, when $\Gamma$ is not torsion-free and say that $\E$, respectively  $\E_0$, is the sheaf associated to the local system on $X$, respectively $X_0$, defined by $(E,\rho)$.

\begin{prop}\Label{prop:pushforward}
The sheaf $\E_0$ on $X_0$ is isomorphic to the push-forward $r_*\E$ of the sheaf $\E$ on $X$.
\end{prop}
\begin{proof}
Let $V_0 \subset X_0$ be a contractible open set.  Let $U_0$ be a connected component of $\tilde{V}_0=\pi_0^{-1}(V_0)$.  Since $\Gamma$ acts properly discontinuously on $D_0$, by shrinking $V_0$ (and hence shrinking $U_0$) if necessary, one can arrange that $\Gamma_{U_0}\equiv \{\gamma \in \Gamma\;  |\;  \gamma \cdot U_0 \cap U_0 \neq \emptyset\}$ is a finite group.   A section $f \in \E_0(V_0)$ is a locally constant map $f: \tilde{V}_0 \to E$ such that $f(\gamma \cdot x)=\rho(\gamma)f(x)$ for all $\gamma \in \Gamma$.  Since $V_0$ is a connected open set, it follows that $\pi_0(U_0)=V_0$ and $\Gamma \cdot U_0 =\pi^{-1}_0(V_0)$, and hence $f$ is determined by its value $v_{U_0}$ on $U_0$.  Furthermore, by the $\Gamma$-equivariance of $f$, $v_{U_0}\in E^{\Gamma_{U_0}}$, the subspace of $E$ fixed by $\Gamma_{U_0}$.  

Recall that $r_*\E(V_0)=\E(r^{-1}(V_0))$.  Since $\pi_0 \circ \tilde{r} =r \circ \pi$, there is a connected component $U$ of $\pi^{-1}(r^{-1}(V_0))$ such that $U=\tilde{r}^{-1}(U_0)$.  Since $\Gamma$ acts properly discontinuously on $D$, $\Gamma_U \equiv \{\gamma \in \Gamma\;  |\;  \gamma \cdot U \cap U \neq \emptyset\}$ is finite.  Thus a section $\psi \in r_*\E(V_0)$ is determined by its value $u_U$ on $U$.  Furthermore, by the $\Gamma$-equivariance of $\psi$, $u_U \in E^{\Gamma_U}$, the subspace fixed by $\Gamma_U$.  

Thus to show that the sheaves are isomorphic, it suffices to show that for sufficiently small $V_0\subset X_0$, the groups $\Gamma_{U_0}$ and $\Gamma_U$ defined above are equal.  It is clear that $\Gamma_{U_0} \subseteq \Gamma_U$.  To show the opposite inclusion, suppose $\gamma \in \Gamma_U$.  Then by definition $\gamma \cdot U \cap U \neq \emptyset$.  Let $x\in \gamma \cdot U \cap U$.  In particular, $x\in U$, so $\tilde{r}(x)\in U_0$.  Furthermore, since $x\in \gamma \cdot U$, $\gamma^{-1}\cdot x \in U$.  Thus $\tilde{r}(\gamma^{-1}\cdot x)\in U$ and the $\Gamma$-equivariance of $\tilde{r}$ implies that $\tilde{r}(x)\in \gamma \cdot U_0$ and hence $\gamma \in \Gamma_{U_0}$.
\end{proof}
\end{subsection}

\begin{subsection}{Cohomology of subspaces}
We recall without proof two results of sheaf cohomology.  A reference for this section is \cite{Bre}.
\begin{thm}[{\cite[Theorem 10.6]{Bre}}]\Label{thm:nbd}
Let $\F$ be a sheaf on a paracompact space $X$, and let $A \subset X$ be a closed subspace.  Let $\N$ be the set of all open neighborhoods of $A$.  Then
\[\varinjlim_{N\in \N} H^*(N;\F|_N)\cong H^*(A;\F|A).\]
\end{thm}

\begin{thm}[{\cite[Theorem 11.7]{Bre}}]\Label{thm:vietoris}
Let $X$ be a paracompact space.  Let $f:X \to Y$ a closed map and $\F$ a sheaf on $X$.  Suppose that $H^p(f^{-1}(y);\F|_{f^{-1}(y)})=0$ for $p >0$ and all $y \in Y$.  Then the natural map
\[f^{\dagger}:H^*(Y;f_* \F) \to H^*(X;\F)\] 
is an isomorphism.
\end{thm}

Let $r:X \to X_0$ denote the deformation retraction arising from the $\Gamma$-equivariant deformation retraction of the symmetric space $D$ to the spine $D_0$.  Let $(E,\rho)$ be a $\Gamma$-module and let $\E$ and $\E_0$ denote the associated local systems on $X$ and $X_0$ respectively.  
\begin{thm}\label{thm:retcoh}
$H^*(X;\E) \cong H^*(X_0;\E_0)$.
\end{thm}

\begin{lem}\Label{lem:stabfiber}
For every $z\in D$, 
\begin{align*}
\Stab{z}&= \Stab{\tilde{r}_t(z)}\quad \text{for $t<1$ and}\\
\Stab{z}&\subseteq \Stab{\tilde{r}_1(z)}.  
\end{align*}
\end{lem}

\begin{proof}
Let $\gamma$ be an element of $\Stab{z}$.  Then $\gamma \cdot \tilde{r}_t(z) =\tilde r_t(\gamma \cdot z)=\tilde r_t(z)$.  Hence $\gamma\in \Stab{\tilde{r}_t(z)}$.  Notice that for each $z\notin D_0$, $c(t)=\tilde r_t(z), 0 \leq t \leq 1$ is a reparameterization of a geodesic, and $\Gamma$ acts by isometries.  Thus every $\gamma \in \Stab{z}$ fixes the geodesic through $\tilde r_1(z)$ and $z$.  In particular, $\Stab{z'}=\Stab{z}$ whenever $z=\tilde r_t(z')$ for some $t<1$.
\end{proof}

\begin{proof}[Proof of Theorem~\ref{thm:retcoh}]
We will show that $H^*(N;\E|_N)\cong H^*(r^{-1}(y);\E|_{r^{-1}(y)})$ for all contractible neighborhoods $N$ of $y \in X_0$.  Let $U\subset D$ be a connected component of $\pi^{-1}(r^{-1}(y))$.  Let $\tilde{y}$ denote the unique lift of $y$ in $U$.  Note that $U$ is a finite union of geodesic rays emanating from $\tilde{y}$.  By Proposition~\ref{lem:stabfiber}, the stabilizer groups $\Stab{\tilde x}$ are isomorphic for $\tilde x$ in a connected component of $U \setminus \tilde{y}$.  Let $\tilde{N} \subseteq U$ be a contractible open subset containing $\tilde{y}$ and let $N=\pi(\tilde N)$.  Then it is clear that $H^*(r^{-1}(y);\E|_{r^{-1}(y)}) \cong H^*(N,\E|_N)$.  Apply Theorem~\ref{thm:nbd} with $X=r^{-1}(y),\F=\E|_{r^{-1}(y)}, A=\{y\}$, and $\N$ the family of all contractible open sets containing $y$.  Then 
\[H^*(r^{-1}(y);\E|_{r^{-1}(y)}) \cong \varinjlim_{N\in \N} H^*(N;\E|_{N}) \cong H^*(y;\E|_y).\]
In particular, $H^p(r^{-1}(y);\E|_{r^{-1}(y)})=0$ for $p>0$.

By Theorem~\ref{thm:vietoris}, this implies that $H^*(X;\E)\cong H^*(X_0;r_*\E)$.  By Proposition~\ref{prop:pushforward}, $r_*\E \cong \E_0$ and the result follows.
\end{proof}
\end{subsection}
\begin{subsection}{Computing the cohomology from the cell structure}
The cellular structure allows us to compute the cohomology $H^*(\Gamma\backslash D_0;\E)$ combinatorially once $\Gamma$-representatives of cells and the stabilizers of those cells have been computed.  We first set up some notation.

For each $p$, fix a set $R_p$ of representatives of $\Gamma$-conjugacy classes of  $p$-cells of $D_0$.  Let $[\tau]$ denote the representative of the class of $\tau$.  For each representative $[\tau]$, fix a distinguished maximal flag of cells 
\[F_{[\tau]}=\tau_0<\tau_1<\cdots<\tau_p,\quad \text{where $\tau_i$ is an $i$-cell and $\tau_p=[\tau]$.}\]
For each cell $\tau$, let $\gamma_\tau \in \Gamma$ be such that $\gamma_\tau \cdot F_{[\tau]}$ terminates in $\tau$.  Note that $\gamma_\tau$ is well-defined up to $\Stab{\tau}$.  

For each fixed $p$-cell, let $S_\sigma$ denote the simplicial complex arising from the poset of cells in $\sigma$ with the partial ordering derived from containment.  In particular, the vertices of $S_\sigma$ are the cells contained in $\sigma$ and the $k$-simplices are the $(k+1)$-flags $\sigma_0<\sigma_1<\cdots <\sigma_k$.  Define a map $n_\sigma: \text{\{$p$-simplices of $S_\sigma$\}}\to \{\pm 1\}$ by the equation $\partial S_\sigma=\sum_{F\in S_\sigma} n_\sigma(F)\partial F$.  Multiply by $-1$ if necessary so that $n_{\sigma}(\gamma_\sigma F_{[\sigma]})=(-1)^p$.  Then for each $\sigma$, define a map
\begin{equation}
\label{eq:sgn} 
\sgn_\sigma(\tau)=n_\sigma(\gamma_\tau F_{[\tau]}<\sigma).\end{equation}

\begin{thm}\label{thm:cohomology}
The cohomology $H^*(X_0;\E)$ can be computed from the complex
\[0 \to \bigoplus_{\sigma \in R_0} E^{\Stab{\sigma}} \to \bigoplus_{\sigma \in R_1} E^{\Stab{\sigma}} \to \cdots \to \bigoplus_{\sigma \in R} E^{\Stab{\sigma}}\to 0,\]
where the differential 
\[d:\bigoplus_{\sigma \in R_{p-1}} E^{\Stab{\sigma}} \to \bigoplus_{\sigma \in R_{p}} E^{\Stab{\sigma}}\]
is given by
\[\left(dv\right)_\sigma=\sum_{\tau \text{ a $(p-1)$-cell }\in \partial \sigma} \sgn_\sigma(\tau) \rho\left(\gamma_\tau \right) v_{[\tau]}.\]
Here $\gamma_\tau$, $\sgn$, and $[\cdot]$ are defined above and the vector $v_{[\tau]}$ is the $[\tau]$-component of the vector $v\in \bigoplus_{\sigma \in R_{p-1}} E^{\Stab{\sigma}}$.  
\end{thm}
\end{subsection}
\end{section}
\begin{section}{Cohomology}\Label{sec:cohcomputation}
Recall that $\epsilon,~w,~\sigma,~\tau,$ and $\xi$ were explicitly defined elements of $G$ given in Section~\ref{sec:su21background}.
\begin{thm} \Label{thm:SU(2,1)cohomology}
Let $E$ be a $\Gamma$-module with the action of $\Gamma$ given by $\rho: \Gamma \to \GL(E)$.  Then $H^*(\Gamma \backslash D; \E)$ can be computed from the following cochain complex.  
\[0 \to C^0 \to C^1 \to C^2 \to C^3 \to 0,\]
where 
\begin{align*}
C^0&=E^{\rho(\tau \epsilon w)} \oplus E^{\langle \rho(\epsilon), \rho(w)\rangle} \oplus E^{\rho(\epsilon w)} \oplus E^{\langle \rho(\epsilon w),\rho(\xi^2)\rangle} \oplus E^{\langle \rho(\epsilon w) , \rho(\sigma \epsilon^2)\rangle} \oplus E^{\rho(\xi)}\\
C^1&=E \oplus E^{\rho(\sigma \epsilon w \sigma^{-1})} \oplus E^{\rho(\epsilon w)} \oplus E \oplus  E^{ \rho(\epsilon w)} \oplus E^{ \rho(\epsilon w)} \oplus E^{\rho(\epsilon) } \oplus E^{ \rho(\xi^2) } \oplus E\\
C^2&=E \oplus E \oplus E \oplus  E^{ \rho(\epsilon w)} \oplus E \oplus E \oplus E \\
C^3&=E \oplus E.
\end{align*}
Then for $(\kappa_i) \in C^0$, $(\lambda_i) \in C^1$, and $(\mu_i) \in C^2$, the differentials are given by
\begin{align*}
d_0(\kappa)&=\begin{pmatrix} -\kappa_1+\rho(\xi^2)\kappa_3\\\rho(\xi) \kappa_3-\kappa_5\\\kappa_3-\kappa_5\\-\kappa_3+\rho(\xi)\kappa_3\\\kappa_3-\kappa_4\\-\kappa_2+\kappa_4\\\kappa_1-\kappa_2\\ \kappa_4-\kappa_6\\ \kappa_3-\kappa_6 \end{pmatrix}\\
d_1 (\lambda)&=\begin{pmatrix}-\lambda_1+\rho(\tau\epsilon w)\lambda_1+\rho(\xi)\lambda_4\\\lambda_2-\lambda_3-\lambda_4\\\lambda_4+\rho(\xi)\lambda_4+\lambda_5-\rho(\xi^2)\lambda_5\\-\rho(\sigma \epsilon^2)\lambda_2+\lambda_3-\lambda_5+\rho(\tau \sigma w \tau \sigma^{-1})\lambda_5-\lambda_6+\rho(\epsilon)\lambda_6\\ \lambda_1-\rho(\xi^2)\lambda_5-\lambda_6+\lambda_7\\-\lambda_5-\lambda_8+\lambda_9\\\lambda_4+\lambda_9-\rho(\xi)\lambda_9\end{pmatrix} \\
d_2(\mu)&=\begin{pmatrix} A\mu_1+B\mu_2-\mu_3-\rho(\tau \sigma \epsilon w^{-1})\mu_3-\mu_4+\mu_5-\rho(\epsilon)\mu_5\\ -\mu_3-\mu_6+\rho(\xi^2)\mu_6+\mu_7+\rho(\xi)\mu_7\end{pmatrix}\end{align*}
where \[
A=\rho(I)+\rho(\tau \epsilon w)+\rho(\tau \epsilon w)^2 
\quad \text {and} \quad B= -\rho(I)+\rho(\sigma \epsilon w\sigma^{-1}) -\rho(\sigma \epsilon w \sigma^{-1} \epsilon w). \]
\end{thm}

An application of Theorem~\ref{thm:SU(2,1)cohomology} is the following.

\begin{cor}\label{cor:generators}
The group $\Gamma=\SU(2,1;\Z[i])$ is generated by $\{\epsilon,w,\tau,\sigma\}$.
\end{cor}

\begin{lem}\label{lem:fix}
Let $H$ be a subgroup of a group $G$.  If $H\neq G$, then there exists a representation $(E,\rho)$ of $G$ such that $E^{\rho(H)} \neq E^{\rho(G)}$.
\end{lem}
\begin{proof}
Consider the representation $(E,\rho)$ of functions $\phi:G/H \to \C$ with the left regular action of $G$ on $\phi$.  The characteristic function $\chi_{eH}$ of the identity coset $eH$ is fixed by $H$, but not fixed by $G$ for $G \neq H$.
\end{proof}

\begin{proof}[Proof of Corollary~\ref{cor:generators}]
Let $\rho:\Gamma \to \GL(E)$ be a representation of $\Gamma$.  Then the cohomology $H^0(X;\E)$ is equal to the global sections $\E(X) \cong E^\Gamma$.  Since $H^0(X;\E)$ is the kernel of $d_0$ given in Theorem~\ref{thm:SU(2,1)cohomology}, $H^0(X;\E)\cong E^{\langle \epsilon,w,\tau,\sigma \rangle} $.  It is well-known that $H^0(X;\E)\cong E^{\rho(\Gamma)}$.  The result then follows from Lemma~\ref{lem:fix}.
\end{proof}

The following corollary follows immediately from the form of $d_2$ given in Theorem~\ref{thm:SU(2,1)cohomology}.

\begin{cor}
Let $E$ be a $\Gamma$-module with action given by $\rho:\Gamma \to \GL(E)$.  Let $\E$ denote the associated sheaf.  Then 
\[\rank(H^3(\Gamma \backslash D;\E)) \leq \rank(E)-\rank(E^{\rho(\epsilon w)}).\]
\end{cor}

\begin{cor}
Let $E$ be a finite dimensional complex representation of $\SU(2,1)$.  Let $E_{ij}$ denote the $i\omega_1+j\omega_2$ weight space of $E$, where $\omega_1$ and $\omega_2$ are the fundamental \textup{(}complex\textup{)} weights of $\mathfrak{sl}_3\C$, the complexification of $\mathfrak{su}(2,1)$.  Then 
\begin{equation*}
\dim H^3(\Gamma\backslash D;\E)\leq \sum_{\substack{i>j\\i\equiv j(4)}}\dim E_{ij}+\sum \dim \left(E_{ii} \cap \ker(I+\rho(w))\right).
\end{equation*}
\end{cor}

\begin{proof}
The dimension of $H^3(\Gamma \backslash D; \E)=2 \dim (E)-\rank(d_2)$.  From the form of $d_2$ given in Theorem~\ref{thm:SU(2,1)cohomology}, it follows that
\begin{equation}\Label{eq:coh}
\dim H^3(\Gamma \backslash D;\E) \leq \dim(E)-\rank(\Phi),
\end{equation} 
where $\Phi:E^{\rho(\epsilon w)} \oplus E \to E$ is the linear map given by $\Phi(v)=v_1+v_2-\rho(\epsilon)v_2$.
Let $E_+$ denote the $(+1)$-eigenspace of $\rho(\epsilon)$ and $E_-$ denote the $(-1)$-eigenspace of $\rho(w)$.  Then \eqref{eq:coh} implies that 
\begin{equation}
\dim H^3(\Gamma \backslash D;\E) \leq \dim(E_+\cap E_-).
\end{equation}
The result follows from translating the eigenvalue condition into a condition on weights of $E$.
\end{proof}

First consider the trivial representation $E=\Z$.  Then $H^*(\Gamma \backslash D;\E)$ is isomorphic to the singular cohomology of $\Gamma \backslash D$.  Proposition~\ref{thm:SU(2,1)cohomology}, allows us to explicitly compute the cohomology. 
\begin{thm}
Let $\Z$ denote the constant sheaf of integers on $\Gamma \backslash D$.  Then \[H^k(\Gamma \backslash D; \Z)=\begin{cases} \Z & k=0,2,\\
                                            0 & k=1 \quad \text{or} \quad  k\geq 3. \end{cases}\]
\end{thm}

\begin{thm}
Let $E=\Z[i]^3$ with the natural action of $\Gamma$ and let $\E$ denote the associated sheaf on $\Gamma \backslash D$.  Then 
\[H^k(\Gamma \backslash D;\E)=\begin{cases} 0 &  k=0,1,\\
\Z^2 & k=2,\\
\Z/2\Z & k=3.\\
\end{cases}\]
\end{thm}

\begin{thm}  The dimensions $h^i$ of the cohomology groups for $\Gamma$ with coefficients in $\Sym^n(V)$, $0 \leq n \leq 20$, where $V\equiv \C^3$ is the standard representation is as tabulated in Table~\ref{tab:coh}.
\begin{table}\Label{tab:coh}
\caption{Cohomology for $\Sym^n(V)$}
\begin{tabular}{|c|cccccccccc |}
\hline
$n$&$1$&$2$&$3$&$4$&$5$&$6$&$7$&$8$&$9$&$10$\\\hline
$h^0$&$0$&$0$&$0$&$0$&$0$&$0$&$0$&$0$&$0$&$0$\\
$h^1$&$0$&$0$&$0$&$0$&$0$&$0$&$0$&$0$&$0$&$0$\\
$h^2$&$1$&$0$&$0$&$1$&$3$&$1$&$2$&$2$&$5$&$1$\\
$h^3$&$0$&$0$&$0$&$1$&$0$&$0$&$0$&$1$&$0$&$0$\\\hline\hline
$n$&$11$&$12$&$13$&$14$&$15$&$16$&$17$&$18$&$19$&$20$\\\hline
$h^0$&$0$&$0$&$0$&$0$&$0$&$0$&$0$&$0$&$0$&$0$\\
$h^1$&$0$&$1$&$0$&$0$&$3$&$4$&$2$&$5$&$8$&$11$\\
$h^2$&$2$&$3$&$7$&$4$&$5$&$4$&$9$&$5$&$7$&$5$\\
$h^3$&$0$&$1$&$0$&$0$&$0$&$1$&$0$&$0$&$0$&$1$\\\hline
\end{tabular}
\end{table}
%
\end{thm}
\end{section}
\bibliography{../references}    
\end{document}